\documentclass{mathincs}
\usepackage[latin1]{inputenc}
\usepackage[english]{babel}
\usepackage{url}
\usepackage{amsfonts,cyracc}
 \newtheorem{thm}{Theorem}[section]
 \newtheorem{cor}[thm]{Corollary}
 \numberwithin{equation}{section}
 \newcommand\R{\mathbb{R}}

\begin{document}
%
%
%
\title[Isochronicity conditions for some planar polynomial systems II ]
 {Isochronicity conditions for some planar polynomial systems II
}
\author{Magali Bardet}
\address{LITIS, Universit\'e de Rouen \\
 Avenue de l'Université BP 12\\
 F-76801 Saint Etienne du Rouvray, \;France}
 \email{magali.bardet@univ-rouen.fr}
 
\author{Islam Boussaada}
\address{Laboratoire de Mathématiques Rapha\"el Salem,\\
 CNRS, Universit\'e de Rouen \\
 Avenue de l'Université BP 12\\
 76801 Saint Etienne du Rouvray, \;France}
\email{islam.boussaada@gmail.com}

\author{A. Raouf Chouikha}
\address{   	
 	Laboratoire Analyse G\'eometrie et Aplications,\\
	 UMR CNRS 7539, Institut Gallilée, Université Paris 13,\\
	  99 Avenue J.-B. Clément, 93430 Villetaneuse, France}
\email{chouikha@math.univ-paris13.fr}

\author{Jean-Marie Strelcyn}
\address{Laboratoire de Mathématiques Rapha\"el Salem,\\
 CNRS, Universit\'e de Rouen \\
 Avenue de l'Université BP 12\\
 76801 Saint Etienne du Rouveray, \;France}
\email{ Jean-Marie.Strelcyn@univ-rouen.fr\medskip}
\address{ Laboratoire Analyse G\'eometrie et Aplications,\\
	 UMR CNRS 7539, Institut Gallilée, Université Paris 13,\\
	  99 Avenue J.-B. Clément, 93430 Villetaneuse, France}
\email{strelcyn@math.univ-paris13.fr}


\maketitle {}

\begin{abstract}
  We study the isochronicity of centers at $O\in \mathbb{R}^2$ for
  systems
  $$\dot x=-y+A(x,y),\;\dot y=x+B(x,y),$$ where $A,\;B\in \mathbb{R}[x,y]$,
  which can be reduced to the Li\'enard type equation. When $deg(A)\leq 4$ and $deg(B) \leq 4$, using the
  so-called C-algorithm we found $36$ new multiparameter families of isochronous
  centers.  For a large class of isochronous centers we provide an explicit general formula for linearization. This paper is a direct continuation of \cite{BoussaadaChouikhaStrelcyn2010} but can be read independently.
\end{abstract}

\section{Introduction}

Let us consider the system of real differential equations of the form
\begin{equation}\label{GEN1}
\frac{dx}{dt}=\dot x=-y+A(x,y),\qquad \frac{dy}{dt}=\dot y=x+B(x,y),
\end{equation}
where $(x,y)$ belongs to an open connected subset $U\subset {\R}^2$
containing the origin $O=(0,0)$, with $A,B\in{C}^1(U,\R)$ such that $A$
and $B$ as well as their first derivatives vanish at $O$.  An isolated
singular point $p\in U$ of system~\eqref{GEN1} is a {\it{center}} if
there exists a punctured neighborhood $V\subset U$ of $p$ such that
every orbit of~\eqref{GEN1} lying in $V$ is a closed orbit surrounding
$p$.  A center $p$ is {\it{isochronous}} if the period is constant for
all closed orbits in some neighborhood of $p$.
  
The simplest example is the linear isochronous center at the origin
$O$ given by the system 
\begin{equation}\label{LINC} 
  \dot x=-y,\; \dot  y=x.
\end{equation} 
The problem of characterization of couples $(A,B)$ such that $O$ is an
isochronous center (even a center) for the system \eqref{GEN1} is
largely open.
  
The well known Poincar\'e Theorem asserts that when $A$ and $B$ are
real analytic, a center of~\eqref{GEN1} at the origin $O$ is
isochronous if and only if in some real analytic coordinate system it
takes the form of the linear center~\eqref{LINC} (see for example
\cite{AmelkinLukashevichSadovskii1982}, Th.13.1, and
\cite{RomanovskiShafer2009}, Th.4.2.1).

An overview \cite{ChavarrigaSabatini1999} presents the basic results
concerning the problem of the isochronicity, see also
\cite{AmelkinLukashevichSadovskii1982,Chouikha2005-1,Chouikha2005-2,Chouikha2007,RomanovskiShafer2009}. 
As this paper is a direct continuation of
\cite{BoussaadaChouikhaStrelcyn2010}, we refer the reader to it for
general introduction to the subject. Here we will recall only the
strictly necessary facts.
 
In some circumstances system \eqref{GEN1} can be reduced to {\it{the
    Li\'enard type equation}}
\begin{equation}\label{L2} 
  \ddot x+f(x){\dot x}^2+g(x)=0
\end{equation} 
with $f,\;g \in C^1(J,\mathbb{R})$, where $J$ is some neighborhood of
$0\in\mathbb{R}$ and $g(0)=0$. In this case, system~\eqref{GEN1} is
called {\it{reducible}}.  To equation~\eqref{L2} is associated the
equivalent, two dimensional (planar), Li\'enard type system
\begin{equation}\label{LienardType}
  \left.\begin{aligned} \dot x &= y \\ \dot y &= -g(x) - f(x) y^2
    \end{aligned}\right\}. 
\end{equation}   
For reducible systems considered in this paper, the nature (center and
isochronicity) of the singular point $O$ for both systems~\eqref{GEN1}
and ~\eqref{LienardType} is the same.
 

Let us return now to the Li\'enard type equation \eqref{L2}. Let us
define the following functions
\begin{equation}\label{xixi}
F(x):= \int_0^xf(s) ds, \quad \phi(x):= \int_0^x e^{F(s)} ds.
\end{equation}
The first integral of the system~\eqref{LienardType} is given by the formula (\cite{Sabatini2004}, Th.1)
 \begin{equation}\label{L2FI}
 I(x,\dot x)=\frac{1}{2}(\dot x e^{F(x)})^2+\int_0^x g(s) e^{2F(s)} ds.
 \end{equation}
When $xg(x)>0$ for $x\neq0$, define the function $X$  by
\begin{equation}\label{xi}
\frac {1}{2} \xi(x)^2 = \int_0^x g(s) e^{2F(s)} ds
\end{equation}
and $x\xi(x)>0$ for $x\neq 0$.

\begin{thm}[\cite{Sabatini2004}, Theorem~2]\label{thm:Center}
  Let $f ,\, g \in C^1(J,\R)$. If $xg(x)>0$ for $x\neq0$, then the
  system~\eqref{LienardType} has a center at the origin $O$.  When
  $f$ and  $g$ are real analytic, this condition is also necessary.
\end{thm}

\begin{thm}[\cite{Chouikha2007}, Theorem~2.1]\label{thm:URABE}
  Let $f$ and $g$ be real analytic functions defined in a neighborhood
  $J$ of $0\in \mathbb{R}$, and let $x g(x) > 0$ for $x \neq 0$.  Then
  system~\eqref{LienardType} has an isochronous center at $O$ if and only if
  there exists an odd function $h\in C^1(J,\R)$ which satisfies the
  following conditions
  \begin{equation}
    \label{CRI}\frac {\xi(x)}{1+h(\xi(x))} = g(x) e^{F(x)},
  \end{equation}
  the function $\phi(x)$ satisfying
  \begin{equation}
    \label{bb} \phi(x) = \xi(x) + \int_0^{\xi(x)} h(t) dt,
  \end{equation}
  and $\xi(x)\phi (x) > 0$ for $x\neq 0$.

  When those equivalent conditions are satisfied, then the function
  $h$ is analytic in the neighborhood $J$. It is called the
  {\emph{Urabe function}} of system~\eqref{LienardType}.
\end{thm}
Taking into account \eqref{xi}, it is easy to see that \eqref{CRI} and
\eqref{bb} are equivalent.


\begin{cor}[\cite{Chouikha2007}, Corollary 2.4] \label{cor:urabenul}
Let $f$ and $g$ be real analytic
  functions defined in a neighborhood of $0\in\R$, and $x g(x) > 0$
  for $x \neq 0$.  The origin $O$ is an isochronous center of
  system~\eqref{LienardType} with Urabe function $h=0$ if and only if
      \begin{equation}\label{Null}
        g'(x)+g(x)f(x)=1
      \end{equation}
  for $x$ in a neighborhood of $0$.
\end{cor}
In the sequel we shall call the Urabe function of the isochronous center of reducible system \eqref{GEN1} the Urabe function of the corresponding Li\'enard type equation.

In \cite{Chouikha2007} the third author described how to use
Theorem~\ref{thm:URABE} to build an algorithm (C-Algorithm, see
Appendix of \cite{BoussaadaChouikhaStrelcyn2010} for more details) to look for isochronous
centers at the origin for reducible system~\eqref{GEN1}, and applied
it to the case where $A$ and $B$ are polynomials of degree $3$.  This
work was continued in \cite{ChouikhaRomanovskiChen2007} and in \cite{BoussaadaChouikhaStrelcyn2010}.


In this paper we apply the so called {\it Rational C-Algorithm} introduced in \cite{BardetBoussaada2010} which is more powerful then the algorithm used in \cite{BoussaadaChouikhaStrelcyn2010}.

The aim of the present paper is to extend these studies to the following real multiparameter family of polynomial system of differential equations :

 \begin{equation}\label{QUAR}
\left.\begin{aligned} \dot x &= -y+a_{{1,1}}xy+a_{{2,0}}{x}^{2}+a_{{2,1}}{x}^{2}y+a_{{3,0}}{x}^{3}+a_{{3,1}}{x}^{3}y+a_{{4,0}}{x}^{4}\\ \dot y &=x+b_{{0,2}}{y}^{2}+b_{{1,1}}xy+
b_{{2,0}}{x}^{2}+b_{{1,2}}x{y}^{2}+b_{{2,1}}{x}^{2}y+b_{{3,0}}{x}^{3}+b_{{2,2}}x^2{y}^{2}+b_{{3,1}}{x}^{3}y+b_{{4,0}}{x}^{4}
 \end{aligned}\right\}. 
   \end{equation}
   The reported results, which are obtained by Maple computations are reproduced without almost any change to avoid misprints.
   
   The paper is organized as follows. In Section \ref{sec::prem} we report the necessary background and describe the investigated subfamilies of system \eqref{QUAR}. In Sections \ref{sec::strd} and \ref{sec::CG} we describe the obtained new isochronous centers. In total we provide $36$ new families of isochronous centers. Among them two {\it Monsters} \eqref{QUARUN42} and \eqref{QUARUN57} of extreme complexity, never encountered before. 
  
   Let us stress that when describing the Urabe functions of the isochronous centers from Section 3, for the first time we encounter the {\it non-standard} examples of it. Indeed, up to now all identified Urabe functions was always of the form $h(\xi)=\frac{a\,\xi^{2n+1}}{\sqrt{b+c\,\xi^{4n+2}}}$ where $a,\,b,\, c\in \mathbb{R}$, $b>0$ and $n$ a non negative integer.(See \cite{Chouikha2007,ChouikhaRomanovskiChen2007,BoussaadaChouikhaStrelcyn2010})
   
   Finally, in Section \ref{sec::lin}, when Urabe function $h=0$, we describe the explicit  general formula for linearizing change of coordinates whose existence is insured by the Poincar\'e theorem. We report also $5$ examples of such linearization.
   \section{Preliminaries}\label{sec::prem}
   \subsection{Choudhury-Guha Reduction}
   
   Let us consider the  real polynomial  system
   \begin{equation}\label{CHERKAS}
\left.\begin{aligned} \dot x &= p_{{0}} \left( x \right) +p_{{1}} \left( x \right) y \\ \dot y &= q_{{0}} \left( x \right) +q_{{1}} \left( x \right) y+q_{{2}} \left( x
 \right) {y}^{2}
\end{aligned}\right\}, 
  \end{equation} 
   where $p_0,p_1,q_0,q_1,q_2 \in \mathbb{R}[x]$.
   
  We will always assume that $O=(0,0)\in \mathbb{R}^2$ is a singular point of \eqref{CHERKAS}, that is $p_0(0)=q_0(0)=0$.  Let us assume also that $p_1(0)\neq 0$. 
  
  Let us note that the system \eqref{QUAR} is a particular case of \eqref{CHERKAS} when
   
\begin{equation}\label{POLYS4}
\left.\begin{aligned}  p_{{0}} \left( x \right) &=a_{{2,0}}{x}^{2}+a_{{3,0}}{x}^{3}
+a_{{4,0}}{x}^{4}\\
p_{{1}} \left( x \right) &=-1+a_{{1,1}}x+a_{{2,1}}{x}^{2}+a_{{3,1}}{x}^{3}\\
q_{{0}}
 \left( x \right) &=x+b_{{2,0}}{x}^{2}+b_{{3,0}}{x}^{3}+b_{{4,0}}{x}^{4}\\
 q_{{1}} \left( 
x \right) &=b_{{1,1}}x+b_{{2,1}}{x}^{2}+b_{{3,1}}{x}^{3}\\
q_{{2}} \left( x \right)& =b_{{0,2}}+b_{{1,2}}x+b_{{2,2}}x^2
 \end{aligned}\right\}. 
   \end{equation} 

The following change of coordinates $x=x,\;z= p_{{0}} \left( x \right) +p_{{1}} \left( x \right) y$ transforms the system \eqref{CHERKAS} to the system
\begin{equation}\label{L2P}
\left.\begin{aligned} \dot x &= z \\ \dot z &= \left( {\frac {q_{{2}} \left( x \right) }{p_{{1}} \left( x \right) }}
+{\frac {p'_{{1}} \left( x \right) }{p_{{1}} \left( x
 \right) }} \right) {z}^{2}+ \left( -{\frac { \left( p'_
{{1}} \left( x \right)  \right) p_{{0}} \left( x \right) }{p_{{1}}
 \left( x \right) }}+q_{{1}} \left( x \right) +p'_{{0}}
 \left( x \right) -2\,{\frac {q_{{2}} \left( x \right) p_{{0}} \left( 
x \right) }{p_{{1}} \left( x \right) }} \right) z\\
&+{\frac {q_{{2}}
 \left( x \right)  \left( p_{{0}} \left( x \right)  \right) ^{2}}{p_{{
1}} \left( x \right) }}-q_{{1}} \left( x \right) p_{{0}} \left( x
 \right) +p_{{1}} \left( x \right) q_{{0}} \left( x \right) 
 \end{aligned}\right\}. 
  \end{equation} 
  If  
\begin{equation}\label{RES}
  -{\frac { \left( p'_
{{1}} \left( x \right)  \right) p_{{0}} \left( x \right) }{p_{{1}}
 \left( x \right) }}+q_{{1}} \left( x \right) +p'_{{0}}
 \left( x \right) -2\,{\frac {q_{{2}} \left( x \right) p_{{0}} \left( 
x \right) }{p_{{1}} \left( x \right) }}=0,\end{equation}
 the system~\eqref{CHERKAS} is of Li\'enard type~\eqref{LienardType},
 
with
 \begin{equation}\label{fg}
\left.\begin{aligned} f(x)&= -\left( {\frac {q_{{2}} \left( x \right) }{p_{{1}} \left( x \right) }}
+{\frac {p'_{{1}} \left( x \right) }{p_{{1}} \left( x
 \right) }} \right) \\ 
 g(x)&= -{\frac {q_{{2}}
 \left( x \right)  \left( p_{{0}} \left( x \right)  \right) ^{2}}{p_{{
1}} \left( x \right) }}+q_{{1}} \left( x \right) p_{{0}} \left( x
 \right) -p_{{1}} \left( x \right) q_{{0}} \left( x \right)
 \end{aligned}\right\} .
\end{equation} 
To the best of our knowledge, the above reduction of the system \eqref{CHERKAS} to Li\'enard type system \eqref{LienardType}  was proposed for the first time in a preliminary and never published version of \cite{ChoudhuryGuha2010}. It is then natural to name it Choudhury-Guha reduction.

In all considered cases (see \eqref{POLYS4}) it is easy to see that
for $|x|$ small enough $g(x)=x+x^2\,\tilde{g}(x)$ where $\tilde{g}$ is
a real analytic function. Thus $x\,g(x)>0$ for $x\neq 0$, $|x|$ small
enough and Theorem~\ref{thm:URABE} insures that the origin $O$ is a
center for the system \eqref{LienardType}.
Our aim is to decide when this center is isochronous.

When $p_0$ and $q_1$ identically vanish, \eqref{RES} is satisfied and we recover the {\it{standard reduction}} from \cite{BoussaadaChouikhaStrelcyn2010} (see Case 1 from Sec. 1 of \cite{BoussaadaChouikhaStrelcyn2010}).
Many particular cases of system~\eqref{QUAR} where studied, using this
standard reduction.   
\begin{enumerate}
   \item In \cite{Chouikha2007} : 
$a_{{1,1}}=a_{{2,0}}=a_{{3,0}}=a_{{3,1}}=a_{{4,0}}=b_{{2,2}}=b_{{4,0}}=b_{{3,1}}=b_{{2,1}}=b_{{1,1}}=0$.
  That means that $p_0(x)=q_1(x)=0$, $p_1(x)=-1+a_{2,1}x^2$ and
  $\deg(q_0(x))\le 3$, $\deg(q_2(x))\le1$.
  \item In \cite{ChouikhaRomanovskiChen2007} : $a_{{2,0}}=a_{{3,0}}=a_{{3,1}}=
   a_{{4,0}}=b_{{2,2}}=b_{{4,0}}=b_{{3,1}}=b_{{2,1}}=b_{{1,1}}=0$. That
   means that $p_0(x)=q_1(x)=0$ and we consider only cubic systems.   
    \item In \cite{BoussaadaChouikhaStrelcyn2010} three families are studied :
    \begin{enumerate}
   \item $a_{{2,0}}=a_{{3,0}}=a_{{4,0}}=b_{{3,1}}=b_{{2,1}}=b_{{1,1}}=0$ with zero Urabe function.  
   \item $a_{{1,1}}=b_{{3,0}}=a_{{2,0}}=a_{{3,0}}=a_{{4,0}}=b_{{3,1}}=b_{{2,1}}=b_{{1,1}}=0$.
    \item $a_{{1,1}}=a_{{2,1}}=a_{{2,0}}=a_{{3,0}}=a_{{4,0}}=b_{{3,1}}=b_{{2,1}}=b_{{1,1}}=0$.
   \end{enumerate}
\end{enumerate}

\subsection{Investigated families}
The exhaustive study of all isochronous center at the origin for the system \eqref{QUAR} is hopeless at the present. Even for cubic system when all quartic terms vanish this problem is not yet solved.

Let us note that the condition \eqref{RES} is equivalent to the
following system of equations:

\begin{equation}\label{RESQUAR}
\left.\begin{aligned}
& 2\,a_{{2,0}}+b_{{1,1}}=0,\\
&-b_{{2,1}}+a_{{1,1}}b_{{1,1}}+a_{{1,1}}a_{{2,0
}}-2\,b_{{0,2}}a_{{2,0}}-3\,a_{{3,0}}=0,\\
&a_{{1,1}}b_{{2,1}}-4\,a_{{4,0}}-
b_{{3,1}}+2\,a_{{1,1}}a_{{3,0}}-2\,b_{{0,2}}a_{{3,0}}-2\,b_{{1,2}}a_{{
2,0}}+a_{{2,1}}b_{{1,1}}=0,\\
&-a_{{3,1}}a_{{2,0}}+3\,a_{{1,1}}a_{{4,0}}+a_{
{2,1}}b_{{2,1}}-2\,b_{{0,2}}a_{{4,0}}-2\,b_{{2,2}}a_{{2,0}}+a_{{2,1}}a
_{{3,0}}+a_{{1,1}}b_{{3,1}}+\\
&
-2\,b_{{1,2}}a_{{3,0}}+a_{{3,1}}b_{{1,1}}=0,\\
&a
_{{3,1}}b_{{2,1}}+a_{{2,1}}b_{{3,1}}+2\,a_{{2,1}}a_{{4,0}}-2\,b_{{1,2}
}a_{{4,0}}-2\,b_{{2,2}}a_{{3,0}}=0,\\
&a_{{3,1}}b_{{3,1}}-2\,b_{{2,2}}a_{{4,0}}+a_{{3,1}}a_{{4,0}}=0.
 \end{aligned}\right\}.
   \end{equation}

In the present paper we determine all isochronous centers of the system \eqref{QUAR} in each of the following three cases.
\begin{enumerate}
\item When the standard reduction is possible (i.e. $p_0(x)=0$ and $q_1(x)=0$, that is
  $a_{{2,0}}=a_{{3,0}}=a_{{4,0}}=b_{{3,1}}=b_{{2,1}}=b_{{1,1}}=0$).
  We provide all candidates for isochronous centers in the cases where
  either $a_{{1,1}}=0$, or $b_{{2,0}}=-3\,b_{{0,2}}$. In all the cases
  but one (a subcase of $b_{{2,0}}=-3\,b_{{0,2}}$), we prove the
  isochronicity. The general case is not yet completely explored.
\item When Choudhury-Guha reduction for the cubic case is possible (i.e. conditions \eqref{RESQUAR}
  are satisfied and
  $a_{{3,1}}=a_{{4,0}}=b_{{3,1}}=b_{{2,2}}=b_{{4,0}}=0$
  ).  For this case we
  obtain the exhaustive list of all isochronous centers at the origin.
\item When Choudhury-Guha reduction is possible and  the Urabe function is null. That means that condition \eqref{RESQUAR} and condition \eqref{Null} for $f$ and
  $g$ defined by \eqref{fg} are simultaneously satisfied. In this case we provide  $25$ examples of new isochronous centers and our analysis is not exhaustive.
\end{enumerate}

Moreover, when the Urabe function $h=0$, we give the explicit formulas for linearizing coordinates from Poincar\'e Theorem. We report $5$ examples where such coordinates are explicitly computed. 
\subsection{Time-reversible systems}
The general notion of  time-reversible system of ordinary differential equations  goes  back to \cite{Devaney1976} where the motivations and general discussion can be found. Here we follow Sec. 1 of \cite{ChavarrigaGineGarcia2001} (see also Sec. 3.5 \cite{RomanovskiShafer2009}). 

The planar system \eqref{GEN1} of ordinary differential equation is {\it time-reversible } if there exists at least one straight line passing through the origin which is a symmetry axis of the phase portrait of the system under consideration. By appropriate rotation this straight line is mapped on the x-axis and the phase portrait of the rotated system is invariant with respect to symmetry $(x,y)\rightarrow(x,-y)$ if only one change time $t$ into $-t$.

Note that a system $$\dot x=P(x,y),\;\dot y=Q(x,y),$$
is time-reversible system with respect to x-axis if and only if $P(x,-y)=-P(x,y)$ and $Q(x,-y)=Q(x,y)$.
 When $P$ and $Q$ are polynomials, this means that the variable y appears in all monomials of $P$  in odd power and   in all monomials of $Q$  in even power ($0$ included).

 Consequently, to decide if a polynomial center for system \eqref{GEN1} is time-reversible or not, we consider the rotated system in coordinates $(x_\alpha,y_\alpha)$, where $x_\alpha= x\, cos \alpha-y \,sin \alpha$ and  $y_\alpha = x\, sin \alpha +y \,cos \alpha$ and we examine the parity of the powers of the variable $y_\alpha$ for all angles $\alpha$.
 
 This notion plays an essential role in our topics. Indeed, for system \eqref{GEN1} the origin is either a center or a focus. Thus, if such system is time-reversible the focus case is excluded and the origin is necessarily a center. 

To the best of our knowledge the majority of already known isochronous centers for polynomial system \eqref{GEN1} are time-reversible. For instance, all systems studied in \cite{Urabe1962,ChavarrigaGineGarcia2001,Chouikha2007,ChouikhaRomanovskiChen2007,BoussaadaChouikhaStrelcyn2010,ChenRomanovski2010} are time reversible. Moreover, among $27$ polynomial isochronous centers presented in tables  $3-29$ of \cite{ChavarrigaSabatini1999} only $7$ are not time-reversible; indeed, those from tables  $17$ and $23-28$.
In what concerns the cubic isochronous centers for system \eqref{GEN1} the complete enumeration of those which are time reversible was obtained in \cite{ChenRomanovski2010}; there are $17$ such cases.
In \cite{ChavarrigaSabatini1999} one find $4$ non time-reversible isochronous centers (tables $25-28$) and in the present paper we present three new such cases which are described in Theorem \ref{CUBICS}. In \cite{ChenRomanovskiZhang2008} a complete list of quartic homogeneous time-reversible isochronous centers is provided, there are $9$ such cases.
In the present paper we provide $33$ new cases of quartic (non homogeneous) isochronous centers. Among them  $8$ are time reversible and at least $23$ are not time reversible.
\subsection{Background on Gr\"obner bases}
The use of the Rational C-Algorithm leads to a system of polynomial
equations 
\begin{equation}
  \label{eq:f1fm}
  f_1=0, \ldots,f_m=0
\end{equation}
with $f_i \in \R[x_1,\ldots,x_n]$. To solve
this system
we consider the ideal $\langle f_1,\ldots,f_m\rangle \subset
\R[x_1,\ldots,x_n]$.
For this aim, we use Gr\"obner Bases computations. In this
section, we recall the basic facts about Gr\"obner bases, and refer the
reader to~\cite{CoxLittleOShea2007} for details.

A
\emph{monomial ordering} is a total order on monomials that is
compatible with the product and such that every nonempty set has a
smallest element for the order. The leading term of a polynomial is
the greatest monomial appearing in this polynomial.

A \emph{Gr\"obner basis} of an ideal~$\mathcal{I}$ for a given
monomial ordering is a set~$G$ of generators of~$\mathcal{I}$ such
that the leading terms of $G$ generate the ideal of leading terms of
polynomials in $\mathcal{I}$. A polynomial is \emph{reduced} with
respect to the Gr\"obner basis~$G$ when its leading term is not a
multiple of those of~$G$. The basis is \emph{reduced} if each
element~$g\in G$ is reduced with respect to~$G\setminus\{g\}$.  For a
given monomial ordering, the reduced Gr\"obner basis of a given set of
polynomials exists and is unique, and can be computed using one's
favorite general computer algebra system, like Maple, Magma or
Singular.  The most efficient Gr\"obner basis algorithm is currently
$F_4$~\cite{Faugere1999}, which is implemented in the three above
cited systems. For our computations, we use the FGb implementation of
$F_4$ available in Maple~\cite{FaugereFGb}.

The complexity of a Gr\"obner basis computation is well known to be generically
exponential in the number of variables, and in the worst case doubly
exponential in the number of variables. Moreover, the choice of the
monomial ordering is crucial for time of the computation.

The \emph{grevlex} ordering is the most suited ordering for the
computation of the (reduced) Gr\"obner basis. The monomials are first
ordered by degree, and the order between two monomials of the same
degree~$x_\alpha=x_1^{\alpha_1}\dotsm x_n^{\alpha_n}$
and~$x_\beta=x_1^{\beta_1}\dotsm x_n^{\beta_n}$ is given
by~$x_\alpha\succ x_\beta$ when the last nonzero element of
$(\alpha_1-\beta_1,\dots,\alpha_n-\beta_n)$ is negative. Thus, among
the monomials of degree~$d$, the order is
\[x_1^d\succ x_1^{d-1}x_2\succ x_1^{d-2}x_{2}^2\succ\dots\succ
x_2^d\succ x_1^{d-1}x_3\succ x_1^{d-2}x_2x_3\succ
x_1^{d-2}x_3^2\succ\dots\succ x_n^d.\] 

However, a Gr\"obner basis for the grevlex ordering is not appropriate
for the computation of the solutions of the system~\eqref{eq:f1fm}.
The most suited ordering for this computation is the
\emph{lexicographical} ordering (or \emph{lex} ordering for
short). The monomials are ordered by comparing the exponents of the
variables in lexicographical order. Thus, any monomial containing
$x_1$ is greater than any monomial containing only variables
$x_2,\ldots,x_n$.

Under some hypotheses (radical ideal with a finite number of
solutions, and up to a linear change of coordinates), the Gr\"obner
Basis of an ideal $\langle f_1,\ldots,f_m\rangle$ for the
lexicographical order $x_1> \ldots > x_n$ has the shape 
\begin{equation}\label{shapelemma}
\{x_1-g_1(x_n), x_2-g_2(x_n),\ldots, x_{n-1}-g_{n-1}(x_{n-1}), g_n(x_n)\},
\end{equation}
where the $g_i$'s are univariate polynomials. In this case, the
computation of the solutions of the system follows easily. In the
general case, the shape of the Gr\"obner basis for the lexicographical
ordering is more complicated, but it is equivalent to several
triangular systems for which the computation of the solutions are
easy. 

An important point is that a Gr\"obner basis for the lex order is in
general hard to compute directly. It is much faster to compute first a
Gr\"obner basis for the grevlex order, and then to make a change of ordering to
the lex order.

The precise ordering we use to compute the Gr\"obner bases of the
polynomial systems occuring in this paper is a weighted order: we fix
a weight $i+j-1$ for the variables $a_{i,j}$ and $b_{i,j}$
(see~\eqref{QUAR}), and use the weighted grevlex or lex
ordering. For those orderings, the polynomials are {\tt homogeneous},
which simplifies the computation. Indeed, without loss of generality,
we can pick a variable $a_{i,j}$ and split the computations into two
cases $a_{i,j}=0$ and $a_{i,j}=1$ (the same concerns $b_{i,j}$).  The
entire set of solutions can then be recovered in the standard way.
For instance, all solutions with $a_{1,1}\neq 0$ for
system~\eqref{QUAR} are obtained from solutions with $a_{1,1}=1$ by
the change of variables $X=a_{1,1}\,x$ and $Y=a_{1,1}\,y$.
This trick reduces by one the number of variables for the Gr\"obner
basis computation and improves the time of the computations. {\bf In what
follows, all results are presented up to such homogenization}.

Finally, we use repeatedly the {\it Radical Membership Theorem}:
\begin{thm}[\cite{CoxLittleOShea2007}]
  Let $I=\langle f_1,\ldots,f_s\rangle$ be an ideal of $k[x_1,\ldots,x_n]$, then $f$ belongs to
  $\sqrt{I}$ if and only if $\langle f_1,\ldots,f_s, 1-yf\rangle =
  \langle 1\rangle = k[x_1,\ldots,x_n,y]$.
\end{thm}


\section{The standard reduction}\label{sec::strd}

In this section we are concerned by system \eqref{QUAR} with  $a_{{2,0}}=a_{{3,0}}=a_{{4,0}}=b_{{3,1}}=b_{{2,1}}=b_{{1,1}}=0$ which gives :

\begin{equation}\label{C_4} 
 \left. \begin{aligned}\dot x&= - y+ a_{1,1}xy+ a_{2,1}x^2y + a_{3,1}x^3y\\
 \dot y&= x + b_{2,0}x^2 + b_{3,0}x^3+ b_{0,2}y^2 + b_{1,2}xy^2 + b_{2,2}x^2y^2+ b_{4,0}x^4 \end{aligned}\right\}. 
\end{equation}

Recall that all cases when the origin $O$ is an isochronous center of
the system \eqref{C_4} with zero Urabe function are described in
\cite{BoussaadaChouikhaStrelcyn2010}. In the following theorem we omit
all isochronous centers with zero Urabe function, as well as all cubic
isochronous centers that where all described
in~\cite{ChouikhaRomanovskiChen2007,Chouikha2007}.

For each case we prove the isochronicity by determining explicitly
its Urabe function. For system $\left(k\right)$ we will denote it by
$h_{\left(k\right)}$.

\begin{thm} 

The following particular cases of system~\eqref{C_4} have an
isochronous center at the origin $O$.
\begin{equation}\label{ST11} 
  \left. \begin{aligned}\dot x&= -y+3\,{x}^{2}y\pm \sqrt {2}{x}^{3}y\\
      \dot y&= x\pm\sqrt {2}{x}^{2}\mp \frac{\sqrt {2}}2{y}^{2}+{x}^{3}+4\,x{y}^{2}\pm 2\,\sqrt 
      {2}{x}^{2}{y}^{2}\pm \frac{\sqrt {2}}4{x}^{4}
 \end{aligned}\right\}, \mbox{ where }
h_{\eqref{ST11}}=-{\frac {\pm \xi}{\sqrt {2+9\,{\xi}^{2}}}}.
\end{equation}
\begin{equation}\label{ST13} 
 \left. \begin{aligned}\dot x&= -y+{x}^{3}y\\
 \dot y&= x+\frac12{x}^{2}{y}^{2}-\frac12{x}^{4} \end{aligned}\right\},
\mbox{ where }
h_{\eqref{ST13}}={\frac {{\xi}^{3}}{\sqrt {4+{\xi}^{6}}}}.
\end{equation}
\begin{equation}\label{ST21} 
  \left. 
    \begin{aligned}
      \dot x&=-y-\frac{{x}^{2}y}{2}+\frac{{x}^{3}y}{8}+xy \\
      \dot y&=x-\frac{3\,{x}^{2}}{4}+\frac{{y}^{2}}{4}+{\frac {5\,{x}^{3}}{24}}+\frac{3\,x{y}^{2}}{8}-\frac{{x}^{2}{y}^{2}}{16}-\frac{{x}^{4}}{48}
    \end{aligned}
  \right\}, 
  \mbox{ where } h_{\eqref{ST21}}={\frac {3\,\xi}{\sqrt {16+9\,{\xi}^{2}}}}.
\end{equation}
\begin{equation}\label{ST22} 
  \left. 
    \begin{aligned}
      \dot x&=-y+9\,{x}^{2}y+6\,{x}^{3}y+xy  \\
      \dot y&=x+\frac{3\,{x}^{2}}{2}-\frac{{y}^{2}}{2}+{x}^{3}+12\,x{y}^{2}+12\,{x}^{2}{y}^{2}+\frac{{x}^{4}}{2}
    \end{aligned}
  \right\},
  \mbox{ where } h_{\eqref{ST22}}=-{\frac {\xi}{\sqrt {4+49\,{\xi}^{2}}}}.
\end{equation}
\begin{equation}\label{ST23} 
  \left. 
    \begin{aligned}
      \dot x&=-y- \left(3\,a_{{3,1}}+\frac{2}{9} \right) {x}^{2}y+a_{{3,1}}{x}^{3}y+xy \\
      \dot y&=x+ \left( -3\,a_{{3,1}}+\frac{1}{9} \right) x{y}^{2}
    \end{aligned}
  \right\},
  \mbox{ where } h_{\eqref{ST23}}={\frac {\xi}{\sqrt {(1-27\,a_{{3,1}}){\xi}^{2}+9}}}.
\end{equation}
\begin{equation}\label{ST24} 
  \left. 
    \begin{aligned}
      \dot x&=-y+xy \\
      \dot y&=x-\frac{3\,{x}^{2}}{2}+{y}^{2}+{x}^{3}-\frac{{x}^{4}}{4}
    \end{aligned}
  \right\}, 
  \mbox{ where } h_{\eqref{ST24}}=\frac{\xi}{\sqrt {1+{\xi}^{2}}}.
\end{equation}
 
Moreover, all other possible isochronous centers at $O$ for non cubic
system~\eqref{C_4}, where either $a_{1,1}=1$ or $b_{2,0}=-3\,b_{0,2}$,
and with non vanishing Urabe functions, belong to the family

\begin{equation}\label{ST26} 
 \left. \begin{aligned}
 \dot x&=-y+ \left( -\frac{3}{8}-2\,b_{{2,2}} \right) {x}^{2}y+ \left( \frac{1}{16}+b_{{2,2}}
 \right) {x}^{3}y+xy
  \\
 \dot y&=x-\frac{3\,{x}^{2}}{4}+\frac{{y}^{2}}{4}+\frac{3\,{x}^{3}}{8}-2\,b_{{2,2}}x{y}^{2}+b_{{2,2}
}{x}^{2}{y}^{2}-\frac{{x}^{4}}{16}
 \end{aligned}\right\}. 
\end{equation}

In particular, when $b_{{2,2}} \in \{- \frac{1}{16},\,0,\,
\frac{1}{16}\}$, the origin $O$ is an isochronous center with non-standard Urabe functions:

 \begin{equation*}   
h_{{\{b_{{2,2}}=-\frac{1}{16}\}}}=\frac{\sqrt {2}\sqrt {2\,L
    \left(\frac{{\xi}^{2}}{4} \right) +    8}\sqrt {{\frac
      {{\xi}^{2}}{L \left( \frac{{\xi}^{2}}{4} \right)      }}} \left(
    L \left( \frac{{\xi}^{2}}{4}\right) +3 \right) L  \left(
    \frac{{\xi}^{2}}{4} \right)}{ 2\,{\xi} \left( L \left(
      \frac{{\xi}^{2}}{4} \right) +4 \right) \left( L \left(
      \frac{{\xi}^{2}}{4} \right) +1 \right)} ,
\end{equation*}
where  $L={\it  LambertW}$  is the Lambert function (see \cite{NIST2010}),
\begin{equation*}
h_{{\{b_{{2,2}}=0\}}}=\frac{\sqrt {2}\sqrt {{\frac {-4+{\xi}^{2}+2\,\sqrt {4+2\,{\xi}^{2}}}{{\xi}^{2}}}}\xi\, \left( {\xi}^{2}+2\,\sqrt {4+2\,{\xi}^{2}}+2  \right)}{ \left( 2+{\xi}^{2} \right) \left( \sqrt {4+2\,{\xi}^{2}}+6  \right)} ,
\end{equation*}
\begin{equation*}
h_{{\{b_{{2,2}}=\frac{1}{16}\}}}={\frac {\sqrt {2}\xi\,\sqrt
    {2\,{\xi}^{2}+32} \left( {\xi}^{2}+12    \right) }{ 2\,\left(
      {\xi}^{2}+4 \right) \left( {\xi}^{2}+16    \right) }}.
  \end{equation*} 
 \end{thm}

\begin{proof}
  Necessary conditions.  This part of the proof is based on the
  C-Algorithm. Indeed, $19$ steps are necessary to find the algebraic
  conditions of isochronicity (see Appendix~A
  of~\cite{BoussaadaChouikhaStrelcyn2010}). We did not succeed
  in computing the full grevlex Gr\"obner basis of the corresponding system
  of polynomial equations. We restricted ourselves to the cases
  $a_{1,1}=0$, which gives the cases~\eqref{ST11}-\eqref{ST13}, and
  $\{a_{1,1}=1, b_{2,0}=-3 b_{0,2}\}$ which gives the
  cases~\eqref{ST21}-\eqref{ST23} and~\eqref{ST26}. We also get
  case~\eqref{ST24} as a particular solution.

Sufficient conditions. For each case we determine its Urabe function. For systems \eqref{ST11}-\eqref{ST24} the procedure in Section 2 of \cite{BoussaadaChouikhaStrelcyn2010} is applied.

The search of the Urabe function for sytem \eqref{ST26} is more subtle. Indeed, we verified that for all values of parameters the first $20$ necessary conditions of isochronicity given by C-algorithm are satisfied. This strongly suggests that for all values of parameters the system \eqref{ST26} has an isochronous center at the origin $O$. 
For this system 
$$f(x)=-{\frac {20-96\,xb_{{2,2}}+64\,{x}^{2}b_{{2,2}}-12\,x+3\,{x}^{2}}{
-16-6\,{x}^{2}-32\,{x}^{2}b_{{2,2}}+16\,{x}^{3}b_{{2,2}}+{x}^{3}+
16\,x}}$$ and

 $$g(x)={\frac {1}{256}}\, \left( -16-6\,{x}^{2}-32\,{x}^{2}b_{{2,2}}+16\,{x
}^{3}b_{{2,2}}+{x}^{3}+16\,x \right) x \left( -16+12\,x-6\,{x}^{2}+{
x}^{3} \right). 
$$
from formula \eqref{xi} one obtains
\begin{equation}\label{eq::xicand} 
\xi^2(x) =2\,{x}^{2} \left( 2-x \right) ^{-2\, \left( 16\,b_{{2,2}}+1 \right) 
^{-1}} \left( 4\,{x}^{2}b_{{2,2}}+1/4\,{x}^{2}-x+2 \right) ^{-{
\frac {16\,b_{{2,2}}-1}{16\,b_{{2,2}}+1}}}.
\end{equation}
From formula \eqref{CRI} one deduces that
$$h(\xi(x))=-{\frac {x \left( 12-6\,x+{x}^{2} \right) }{ \left( x-4 \right) 
 \left( {x}^{2}-2\,x+4 \right) }}.$$
Now the problem is to find the reciprocal function $x=x(\xi)$. 
Unfortunately we succeeded in finding it only for $b_{{2,2}} \in \{- \frac{1}{16},\,0,\,\frac{1}{16}\}$ because in those cases the equation \eqref{eq::xicand} takes a  sufficiently simple form.

\end{proof}
Note that the system \eqref{ST13} was already identified in Theorem 2.2 of \cite{BoussaadaChouikhaStrelcyn2010}.
\section{The Choudhury-Guha reduction}\label{sec::CG}
\subsection{Cubic isochronous centers}
Choudhury-Guha reduction is more general than the standard one used in preceding papers \cite{Chouikha2007,ChouikhaRomanovskiChen2007,BoussaadaChouikhaStrelcyn2010}.
Here we provide the complete enumeration of all cubic systems from
\eqref{QUAR} and we find three new cases of isochronous centers at the
origin. The system  we consider is:
\begin{equation}\label{CUB}
\left.\begin{aligned} \dot x &= -y+a_{{1,1}}xy+a_{{2,0}}{x}^{2}+a_{{2,1}}{x}^{2}y+a_{{3,0}}{x}^{3}\\ \dot y &=x+b_{{2,0}}{x}^{2}+b_{{1,1}}xy+b_{{0,2}}{y}^{2}+b_{{2,1}}y{x}^{2}+b_{{
1,2}}{y}^{2}x+b_{{3,0}}{x}^{3}
 \end{aligned}\right\}. 
   \end{equation} 
Condition~\eqref{RES} is equivalent to the following system of equations:
\begin{equation}\label{RESCUB}
\left.\begin{aligned} 
    &2\,a_{{2,0}}+b_{{1,1}}=0\\
 &a_{{1,1}}a_{{2,0}}-3\,a_{{3,0}}-b_{{2,1}}+a_{
{1,1}}b_{{1,1}}-2\,b_{{0,2}}a_{{2,0}}=0\\
 &a_{{1,1}}b_{{2,1}}-2\,b_{{0,2}}a
_{{3,0}}+2\,a_{{1,1}}a_{{3,0}}+a_{{2,1}}b_{{1,1}}-2\,b_{{1,2}}a_{{2,0}
}=0\\
 &a_{{2,1}}a_{{3,0}}+a_{{2,1}}b_{{2,1}}-2\,b_{{1,2}}a_{{3,0}}=0			
 \end{aligned}\right\}. 
   \end{equation}    
    \begin{thm}\label{CUBICS}
      Under the assumptions~\eqref{RESCUB} the origin $O$ is an
      isochronous center of System~\eqref{CUB} only in one of the
      following cases:
  \begin{enumerate}
  \item  The standard reduction is possible, that means $a_{{3,0}}=a_{{2,0}}=b_{{1,1}}=b_{{2,1}}=0$ and the system is one of those from Theorem 3 of~\cite{ChouikhaRomanovskiChen2007}.
\item  We are in one of the following cases :
\begin{equation}\label{CUB1}
\left.\begin{aligned} \dot x &=-y-2\,b_{{2,0}}xy+{x}^{2}+2\,b_{{2,0}}{x}^{3}
\\ \dot y &=x-4\,b_{{2,0}}{y}^{2}-2\,xy+b_{{2,0}}{x}^{2}+4\,b_{{2,0}}{x}^{2}y+2\,{x}^{3}
 \end{aligned}\right\}, 
   \end{equation} 
 
\begin{equation}\label{CUB6} 
\left.\begin{aligned} \dot x &=-y\pm 2\,\sqrt
    {2}xy+{x}^{2}\mp 2\,\sqrt {2}{x}^{3}
\\ \dot y &=x\pm 8\,\sqrt {2}{y}^{2}-2\,xy\mp 3\,\sqrt {2}{x}
^{2}\mp 12\,\sqrt {2}{x}^{2}y+10\,{x}^{3}
 \end{aligned}\right\}, 
   \end{equation} 

\begin{equation}\label{CUB2}
\left.\begin{aligned} \dot x &=-y-\frac{1}{2}\,b_{{2,0}}xy+{x}^{2}+\frac{1}{2}\,b_{{2,0}}{x}^{3}
\\ \dot y &=x-b_{{2,0}}{y}^{2}-2\,xy+b_{{2,0}}{x}^{2}+b_{{2,0}}{x}^{2}y+\left( 2+\frac{1}{4}\,{b_{{2,0}}}^{2} \right) {x}^{3}
 \end{aligned}\right\}. 
   \end{equation} 

   \end{enumerate}
\end{thm} 
\begin{proof}
  The necessary conditions are given by the solutions of the
  polynomial system of equations consisting of
  equations~\eqref{RESCUB} (called $C_1,\ldots,C_4$) and the 8
  equations obtained from the Rational C-Algorithm (15 steps), called
  $C_5,\ldots,C_{12}$. Let us denote by $I$ the ideal generated by
  $C_1,\ldots,C_{12}$.

  We exclude the standard reduction by adding to $I$ the variable $T$
  and the polynomial
   \[C_{13}=(Ta_{3,0}-1)(Ta_{2,0}-1)(Tb_{1,1}-1)(Tb_{2,1}-1).\] For
   $a_{2,0}=0$, a Gr\"obner basis of $\langle
   C_1,\ldots,C_{12},C_{13},a_{2,0}\rangle$ is $\langle 1 \rangle$
   (i.e. there is no solution), which implies that we can take
   $a_{2,0}=1$. We use the weighted order
   $b_{1,1}>b_{2,1}>b_{3,0}>b_{1,2}>a_{2,1}>b_{0,2}>a_{1,1}>b_{2,0}>a_{3,0}$.

   First, a Gr\"obner basis of system $\langle
   C_1,\ldots,C_6\rangle$ for the weighted lex order contains
   the polynomial
   \[P=(a_{1,1}+2\,b_{0,2})(a_{2,1}-a_{1,1}a_{3,0}-a_{3,0}^2).
   \]
   We split our problem into two subcases according to this
   factorization.
   \begin{itemize}
   \item for $a_{1,1}+2\,b_{0,2}=0$, we get only one real solution
     \begin{equation*}\left.
       \begin{aligned}
         \dot x &= -y + x^2\\
         \dot y &= x -2xy + 2x^3
       \end{aligned}\right\},
     \end{equation*}
     which is a particular case of~\eqref{CUB2} with $b_{2,0}=0$.
   \item for $a_{2,1}-a_{1,1}a_{3,0}-a_{3,0}^2=0$, we eliminate the
     solutions that are not real by adding to the system the
     polynomials $P_i\cdot T_i-1$ for each $P_i$ in
\[
\{16\,a_{3, 0}^2+1,4\,a_{3, 0}^2+9, 4\,a_{3, 0}^2+1, a_{3, 0}^2+4,
a_{3, 0}^2+1, a_{3, 0}^2+16, a_{3, 0}^2+9-4\,b_{2, 0}\,a_{3,
  0}+4\,b_{2, 0}^2\}
\]
that have no real solution. Then, the solutions of the resulting
system are those quoted in the theorem.
   \end{itemize}
	
   Sufficiency. For the cases~\eqref{CUB1} and \eqref{CUB6} we have
   $g'(x)+f(x)g(x)=1$. Hence by Corollary~\ref{cor:urabenul} the
   origin is an isochronous center. Moreover we easily check that
   $h_{\eqref{CUB2}}(\xi)=-\frac{1}{2}\,b_{{2,0}}\xi$.
 \end{proof}

\subsection{Quartic isochronous centers}

Our first target was to identify all isochronous centers at the origin with zero Urabe function for the system \eqref{QUAR}
 \begin{equation*}
\left.\begin{aligned} \dot x &= -y+a_{{1,1}}xy+a_{{2,0}}{x}^{2}+a_{{2,1}}{x}^{2}y+a_{{3,0}}{x}^{3}+a_{{3,1}}{x}^{3}y+a_{{4,0}}{x}^{4}\\ \dot y &=x+b_{{0,2}}{y}^{2}+b_{{1,1}}xy+
b_{{2,0}}{x}^{2}+b_{{1,2}}x{y}^{2}+b_{{2,1}}{x}^{2}y+b_{{3,0}}{x}^{3}+b_{{2,2}}x^2{y}^{2}+b_{{3,1}}{x}^{3}y+b_{{4,0}}{x}^{4}
 \end{aligned}\right\} 
   \end{equation*}
 under the condition \eqref{RESQUAR}. That means finding all values of the $15$ parameters for which the equation $g'(x)+f(x)\,g(x)=1$ is satisfied  where $f$ and $g$ are defined by \eqref{fg} (see Corollary \ref{cor:urabenul}). 
 
When the standard reduction is possible, that means $a_{{4,0}}=a_{{3,0}}=a_{{2,0}}=b_{{1,1}}=b_{{2,1}}=b_{{3,1}}=0$ all the $6$ isochronous centers with zero Urabe function  were described in Theorem 3.1 of~\cite{BoussaadaChouikhaStrelcyn2010}. Otherwise when the Choudhury-Guha reduction needs to be applied the problem becomes  substantially more complicated.
 
 Taking in account the great complexity of the problem we did not succeed in solving it completely. Nevertheless, during our investigations we obtained $25$ new isochronous centers for the system \eqref{QUAR}, two of them of extreme complexity, called Monsters. We are convinced that our list is not exhaustive.

 The procedure to obtain the isochronous centers listed below consists in solving  by Gr\"obner method the system  \eqref{RESQUAR} simultaneously with the set of equations on parameters which corresponds the equation  $g'(x)+f(x)\,g(x)=1$. First, one applies the Solve routine of Maple (based on Gr\"obner basis technic) which splits the variety of solutions into $37$ subvarieties. The cases \eqref{QUARUN1}-\eqref{QUARUN42} were obtained by detailed inspection of some of them. The remaining $7$ isochronous centers  \eqref{QUARUN51}-\eqref{QUARUN57} were obtained by restricting ourselves to $b_{2,2}=a_{3,1}=0$ and by application of the standard Gr\"obner basis technique.
 
 We verified also that all above isochronous centers are not time-reversible, except perhaps the two Monsters \eqref{QUARUN42} and \eqref{QUARUN57}.

    \begin{thm}\label{thm:QUARTICURABENUL}
 The following quartic systems have an isochronous center at the origin $O$ with zero Urabe function.

\begin{equation}\label{QUARUN1}
\left.\begin{aligned} \dot x &=-y+b_{{0,2}}xy+{x}^{2}-b_{{0,2}}{x}^{3}
\\                  \dot y &=x+b_{{0,2}}{y}^{2}-2\,xy+2\,{x}^{3}-b_{{0,2}}{x}^{4}
 \end{aligned}\right\}, 
   \end{equation} 
     $$ ~$$ 
     
  \begin{equation}\label{QUARUN12}
\left.\begin{aligned} \dot x &=-y+xy-a_{{3,0}}{x}^{2}+a_{{3,0}}{x}^{3}
\\                  \dot y &=x+{y}^{2}+2\,a_{{3,0}}xy+2\,{a_{{3,0}}}^{2}{x}^{3}-{a_{{3,0}}}^{2}{x}^
{4}
 \end{aligned}\right\}, 
   \end{equation} 
     $$ ~$$ 
     
  \begin{equation}\label{QUARUN14}
\left.\begin{aligned} \dot x &=-y+xy-a_{{3,0}}{x}^{2}+a_{{3,0}}{x}^{3}
\\                  \dot y &=x+3\,{y}^{2}+2\,a_{{3,0}}xy-{x}^{2}+4\,a_{{3,0}}{x}^{2}y+ \left( \frac{1}{3}+2
\,{a_{{3,0}}}^{2} \right) {x}^{3}+{a_{{3,0}}}^{2}{x}^{4}
 \end{aligned}\right\}, 
   \end{equation} 
     $$ ~$$ 
     
  \begin{equation}\label{QUARUN15}
\left.\begin{aligned} \dot x &=-y+xy-a_{{3,0}}{x}^{2}+a_{{3,0}}{x}^{3}
\\                  \dot y &=x+4\,{y}^{2}+2\,a_{{3,0}}xy-\frac{3\,{x}^{2}}{2}+6\,a_{{3,0}}{x}^{2}y+ \left( 
1+2\,{a_{{3,0}}}^{2} \right) {x}^{3}+ \left( -\frac{1}{4}+2\,{a_{{3,0}}}^{2}
 \right) {x}^{4}
 \end{aligned}\right\}, 
   \end{equation}
   $$ ~$$ 
   
  \begin{equation}\label{QUARUN3}
\left.\begin{aligned} \dot x &=-y+\frac{b_{{0,2}}xy}{3}+{x}^{2}-\frac{b_{{0,2}}{x}^{3}}{3}
\\                  \dot y &=x+b_{{0,2}}{y}^{2}-2\,xy-\frac{b_{{0,2}}{x}^{2}}{3}-\frac{4\,b_{{0,2}}{x}^{2}y}{3}+
 \left( \frac{{b_{{0,2}}}^{2}}{27}+2 \right) {x}^{3}+\frac{b_{{0,2}}{x}^{4}}{3}
 \end{aligned}\right\}, 
   \end{equation} 
     $$ ~$$ 
     
  \begin{equation}\label{QUARUN4}
\left.\begin{aligned} \dot x &=-y+\frac{b_{{0,2}}xy}{4}+{x}^{2}-\frac{b_{{0,2}}{x}^{3}}{4}
\\                  \dot y &=x+b_{{0,2}}{y}^{2}-2\,xy-\frac{3\,b_{{0,2}}{x}^{2}}{8}-\frac{3\,b_{{0,2}}{x}^{2}y}{2}+
 \left( \frac{{b_{{0,2}}}^{2}}{16}+2 \right) {x}^{3}+ \left( -{\frac {1}{
256}}\,{b_{{0,2}}}^{3}+\frac{b_{{0,2}}}{2} \right) {x}^{4}
 \end{aligned}\right\}, 
   \end{equation}
     $$ ~$$  
     
  \begin{equation}\label{QUARUN5}
\left.\begin{aligned} \dot x &=-y-{\frac {45}{8}}\,{x}^{2}y+{x}^{2}+{\frac {45}{8}}\,{x}^{4}
\\                  \dot y &=x-2\,xy-{\frac {225}{8}}\,x{y}^{2}+\frac{19}{2}\,{x}^{3}+45\,{x}^{3}y
 \end{aligned}\right\}, 
   \end{equation} 
     $$ ~$$ 
     
  \begin{equation}\label{QUARUN6}
\left.\begin{aligned} \dot x &=-y+\frac{{b_{{4,0}}}^{2}{x}^{2}y}{2}+{x}^{2}-\frac{{b_{{4,0}}}^{2}{x}^{4}}{2}
\\                  \dot y &=x+b_{{4,0}}{y}^{2}-2\,xy-\frac{b_{{4,0}}{x}^{2}}{2}+{b_{{4,0}}}^{2}x{y}^{2}
-2\,b_{{4,0}}{x}^{2}y+2\,{x}^{3}-{b_{{4,0}}}^{2}{x}^{3}y+b_{{4,0}}{x}^
{4}
 \end{aligned}\right\}, 
   \end{equation} 
     $$ ~$$ 
     
  \begin{equation}\label{QUARUN7}
\left.\begin{aligned} \dot x &=-y+ \left( -2+2\,\sqrt {19} \right) {x}^{2}y+{x}^{2}- \left( -2+2\,
\sqrt {19} \right) {x}^{4}
\\                  \dot y &=x\pm \alpha_1{y}^{2}-2\,xy\mp\frac{\alpha_1{x}^{2}}{2}+ \alpha_2 x{y}^{2}\mp2\,\alpha_1{x}^{2}y+ \alpha_4 {x}^{3}+
\alpha_3 {x}^{3}y\pm4\,\alpha_1{x}^{4}
 \end{aligned}\right\}, 
   \end{equation} 
 
   \rm{where}
   
    $\alpha_1=\sqrt {-106+34\,\sqrt {19}},\,\alpha_2= -10+10\,\sqrt {19},\, \alpha_3=  16-16\,\sqrt {19} ,\, \alpha_4= -13+3\,\sqrt {19}$.
      $$ ~$$ 
      
  \begin{equation}\label{QUARUN9}
\left.\begin{aligned} \dot x &=-y+a_{{1,1}}xy+{\frac {15}{8}}\,{a_{{1,1}}}^{2}{x}^{2}y+{x}^{2}-{x}^{3
}a_{{1,1}}-{\frac {15}{8}}\,{x}^{4}{a_{{1,1}}}^{2}
\\                  \dot y &=x-\frac{a_{{1,1}}{y}^{2}}{2}-2\,xy+\frac{3\,a_{{1,1}}{x}^{2}}{4}+{\frac {15}{4}}\,{
a_{{1,1}}}^{2}x{y}^{2}+3\,a_{{1,1}}{x}^{2}y+2\,{x}^{3}-{\frac {15}{4}}
\,{a_{{1,1}}}^{2}{x}^{3}y-\frac{5\,a_{{1,1}}\,{x}^{4}}{2}
 \end{aligned}\right\}, 
\end{equation} 
  $$ ~$$ 
  
  \begin{equation}\label{QUARUN10}
\left.\begin{aligned} \dot x &=-y \mp{\frac {2}{35}}\,\alpha_{{5}}xy+\alpha_{{6}}{x}^{2}y+{x}^{2}\pm{
\frac {2}{35}}\,\alpha_{{5}}{x}^{3}-\alpha_{{6}}{x}^{4}
\\                  \dot y &=x\pm\frac{\alpha_{{5}}{y}^{2}}{35}-2\,xy\mp{\frac {3}{70}}\,\alpha_{{5}}{x}^{2}
+5\,\alpha_{{6}}x{y}^{2}\mp{\frac {6}{35}}\,\alpha_{{5}}{x}^{2}y+\alpha_{{7
}}{x}^{3}-8\,\alpha_{{6}}\, {x}^{3}y\pm{
\frac {38}{35}}\,\alpha_{{5}}{x}^{4}
 \end{aligned}\right\}, 
   \end{equation}

   
  \noindent\rm{where}
   
    $\alpha_{{5}}=\sqrt {-77798+1162\,\sqrt {4691}},\,\alpha_{{6}}=-{\frac {354}{5}}+\frac{6\,\sqrt {4691}}{5},\,\alpha_{{7}}=-{\frac {2183}{35}}+{\frac {27}{35}}\,\sqrt {4691}.$
  $$ ~$$ 
  \begin{equation}\label{QUARUN16}
\left.\begin{aligned} \dot x &=-y+xy-\frac{3\,a_{{3,0}}{x}^{2}}{4}+a_{{3,0}}{x}^{3}-\frac{a_{{3,0}}{x}^{4}}{4}
\\                  \dot y &=x+3\,{y}^{2}+\frac{3\,a_{{3,0}}xy}{2}-{x}^{2}+\frac{9\,a_{{3,0}}{x}^{2}y}{4}+ \left(\frac{ 1}{3}+{\frac {9}{8}}\,{a_{{3,0}}}^{2} \right) {x}^{3}-\frac{3\,a_{{3,0}}{x}^{
3}y}{4}-\frac{3\,{a_{{3,0}}}^{2}{x}^{4}}{8}
 \end{aligned}\right\}, 
   \end{equation} 
     $$ ~$$ 
     
  \begin{equation}\label{QUARUN18}
\left.\begin{aligned} \dot x &=-y+xy\pm\frac{\sqrt {2}}{2}\,{x}^{2}\mp{\frac{2\,\sqrt {2}}{3}\,{x}^{3}}\pm{\frac{\sqrt {2}}{6}\,{x}^{4}}
\\                  \dot y &=x+6\,{y}^{2}\mp\sqrt
{2}xy-\frac{5}{2}\,{x}^{2}\mp \frac{9}{2}\,\sqrt {2}{x}^{2}y+\frac{13}{3}\,{x}
^{3}\pm \frac{3\,\sqrt {2}}{2}{x}^{3}y-\frac{4}{3}\,{x}^{4}
 \end{aligned}\right\}, 
   \end{equation}
     $$ ~$$ 

  \begin{equation}\label{QUARUN32}
\left.\begin{aligned} \dot x &=-y\mp 2\,\sqrt {2}xy+{x}^{2}\pm 2\,\sqrt {2}{x}^{3}
\\                  \dot y &=x\mp 6\,\sqrt {2}{y}^{2}-2\,xy\pm
2\,\sqrt {2}{x}^{2}\pm 8\,\sqrt {2}{x}^{2}y+
\frac{14}{3}\,{x}^{3}\mp 2\,\sqrt {2}{x}^{4}
 \end{aligned}\right\}, 
   \end{equation}
   $$ ~$$ 

	

     
  \begin{equation}\label{QUARUN39}
\left.\begin{aligned} \dot x &=-y+a_{{1,1}}xy+{a_{{1,1}}}^{2}{x}^{2}y+{x}^{2}-a_{{1,1}}{x}^{3}\\
&-{a_{{1,1}}}^{2}{x}^{4}
\\                  \dot y &=x+3\,a_{{1,1}}{y}^{2}-2\,xy-a_{{1,1}}{x}^{2}+2\,{a_{{1,1}}}^{
2}x{y}^{2}\\
&-4\,a_{{1,1}}{x}^{2}y+2\,{x}^{3}-2\,{a_{{1,1}}}^{2}
{x}^{3}y+a_{{1,1}}{x}^{4}
 \end{aligned}\right\}, 
   \end{equation} 
     $$ ~$$ 
     
  \begin{equation}\label{QUARUN40}
\left.\begin{aligned} \dot x &=-y+a_{{1,1}}xy+3\,{a_{{1,1}}}^{2}{x}^{2}y+{x}^{2}-a_{{1,1}}{x}^{3}-3\,{a_{{1,1}}}^{2}{x}^{4}
\\                  \dot y &=x+4\,a_{{1,1}}{y}^{2}-2\,xy-\frac{3}{2}\,a_{{1,1}}{x}^{2}+6\,{a_{{1,1
}}}^{2}x{y}^{2}-6\,a_{{1,1}}{x}^{2}y+2\,{x}^{3}-6\,{a_{{1,1}}
}^{2}{x}^{3}y+2\,a_{{1,1}}{x}^{4}
 \end{aligned}\right\}, 
   \end{equation} 
     $$ ~$$ 
     
  \begin{equation}\label{QUARUN41}
\left.\begin{aligned} \dot x &=-y+a_{{1,1}}xy+ \left( {a_{{1,1}}}^{2}-\frac{3}{2}\,b_{{0,2}}a_{{1,1}}+\frac{1}{2}\,{b_{{0,2}}}^{2} \right) {x}^{2}y\\
&+{x}^{2}-a_{{1,1}}{x}^
{3}- \left( {a_{{1,1}}}^{2}-\frac{3}{2}\,b_{{0,2}}a_{{1,1}}+\frac{1}{2}\,{b_{{0,2}}}^{
2} \right) {x}^{4}
\\                  \dot y &=x+b_{{0,2}}{y}^{2}-2\,xy+ \left( \frac{1}{2}\,a_{{1,1}}-\frac{1}{2}\,b_{{0,2}
} \right) {x}^{2}+ \left( 2\,{a_{{1,1}}}^{2}-3\,b_{{0,2}}a_{{1,1}}+{b_
{{0,2}}}^{2} \right) x{y}^{2}\\
&+ \left( 2\, \left( {a_{{1,1}}}
^{2}-\frac{3}{2}\,b_{{0,2}}a_{{1,1}}+\frac{1}{2}\,{b_{{0,2}}}^{2} \right) -2
\, \left( 2\,{a_{{1,1}}}^{2}-3\,b_{{0,2}}a_{{1,1}}+{b_{{0,2}}}^{2}
 \right)  \right) {x}^{3}y\\
 &+ \left( -2\,a_{{1,1}}+b_{{0,2}}\right) {x}^{4}+ \left( 2\,a_{{1,1}}-2\,b_{{0,2
}} \right) {x}^{2}y+2\,{x}^{3}
 \end{aligned}\right\}. 
   \end{equation} 
   
  $$ ~$$ 
 \begin{equation}\label{QUARUN42}
\left.\begin{split} \dot x &=-y+Txy+M{x}^{2}y+{x}^{2}-T{x}^{3}-M{x}^{4}\\
\dot y&=x+ P {y}^{2}-2\,xy-\frac{P}{2} {x}^{2}+5\,Mx{y}^{2}-2\,P {x}^{2}y+ S{x}^{3}-8\,M{x}^{3}y+ 4\,P {x}^{4}\end{split}\right\}, 
   \end{equation}

\noindent where
 {\small $T={\it RootOf} \left( 8436480\,{{\it \_Z}}^{8}-151490752\,{{\it \_Z}}^
{6}+799019220\,{{\it \_Z}}^{4}-927412425\,{{\it \_Z}}^{2}-1728684180
 \right) $}
has two real values and
 { \small $M={\it RootOf} ( 54647240223827376\,{{\it \_Z}}^{5}+$
 
 $ \left( 
561508456956838398+111131520938630016\,{T}^{2}-25008718903435904\,{T}^
{4}+1701065448975360\,{T}^{6} \right) 
{{\it \_Z}}^{4}+ \left( 
17547249759911040\,{T}^{6}-5911680053099923488+2296193094798110772\,{T
}^{2}-236397331187044016\,{T}^{4} \right) {{\it \_Z}}^{3}+ \left( -
47291962375532160\,{T}^{6}-32577228112841141760\,{T}^{2}+
1824502903342225392\,{T}^{4}-38028219495361710336 \right) {{\it \_Z}}^
{2}+ \left( 54944501282517005964-141937291355478113601\,{T}^{2}+
26240570003047890492\,{T}^{4}-3151542249299751840\,{T}^{6} \right) {
\it \_Z}+464619529558904831820+500960459979496602351\,{T}^{2}+
20184288412145753196\,{T}^{4}-11100198119724668448\,{T}^{6} )$}\\
has three real values.\\

{\small
\noindent$P={\frac {
1277741426613483475261465072750940888801131320332201}{
78462183804600834942501727375359564182536082155306896}}\,{M}^{3}T+{
\frac {387891696411877338314673950037666400855942884617917}{
8718020422733426104722414152817729353615120239478544}}\,{M}^{2}T$

\noindent$-{
\frac {279706937387076388592164504587266774323327561430923}{
242167233964817391797844837578270259822642228874404}}\,MT-{\frac {
742426510189441463550540802769344501623332700250607}{
58846637853450626206876295531519673136902061616480172}}\,{M}^{3}{T}^{3
}$

\noindent$-{\frac {1946020041492661760957343501655785409975201353232}{
14711659463362656551719073882879918284225515404120043}}\,{M}^{3}{T}^{7
}+{\frac {525193640940998707118304158017605885033259029441052}{
220674891950439848275786108243198774263382731061800645}}\,{M}^{3}{T}^{
5}$

\noindent$+{\frac {103442863691365219385975139449969928492687792379916}{
24519432438937760919531789804799863807042525673533405}}\,{M}^{2}{T}^{5
}-{\frac {164131380112569086879618812047151535419534003325851}{
6538515317050069578541810614613297015211340179608908}}\,{M}^{2}{T}^{3}$

\noindent$-{\frac {367033244357164535298438083037260851819189099856}{
1634628829262517394635452653653324253802835044902227}}\,{M}^{2}{T}^{7}
+{\frac {14771752963262353720563034395097445259447260691808}{
220674891950439848275786108243198774263382731061800645}}\,{M}^{4}{T}^{
5}$

\noindent$-{\frac {54666403014280865159557614303552794593967705728}{
14711659463362656551719073882879918284225515404120043}}\,{M}^{4}{T}^{7
}-{\frac {5151178974136622357147584757325383291000410195430}{
14711659463362656551719073882879918284225515404120043}}\,{M}^{4}{T}^{3
}$

\noindent$+{\frac {4690601171400748048835548400017315870720}{
6303942417299603349700937097152217643374801}}\,{T}^{7}-{\frac {
121307988379654758804095121568312112060864}{
18911827251898810049102811291456652930124403}}\,{T}^{5}$

\noindent$-{\frac {
759721576540272607925488426699587086173816}{
6303942417299603349700937097152217643374801}}\,{T}^{3}+{\frac {
3689908965797406058372519939422856344716780525509}{
9807772975575104367812715921919945522817010269413362}}\,T{M}^{4}$

\noindent$+{
\frac {10767883711615521687536020314030171163247140}{
2101314139099867783233645699050739214458267}}\,T-{\frac {
82822172756429241023335015771633450343406027311056}{
544876276420839131545150884551108084600945014967409}}\,M{T}^{5}$

\noindent$+{
\frac {152574332389107076817431693266653145437149456804321}{
181625425473613043848383628183702694866981671655803}}\,M{T}^{3}+{
\frac {1501312340420936913937036728933692151818809426880}{
181625425473613043848383628183702694866981671655803}}\,M{T}^{7}$
}

{\small 
\noindent$S={\frac {425588095896001129751500929401062088207356}{
510619335801267871325775904869329629113358881}}\,{M}^{4}-{\frac {
5767486187745253514611173391155686944173385}{
1021238671602535742651551809738659258226717762}}\,{M}^{3}$

\noindent$+{\frac {
2663265122906318241956827533449698379402528}{
700438046366622594411215233016913071486089}}-{\frac {
59835284352308243659547352890456360623975739}{
226941927022785720589233735497479835161492836}}\,{M}^{2}$

\noindent$+{\frac {
4276254433644147165491110877451438697226467}{
4202628278199735566467291398101478428916534}}\,M+{\frac {
20238619233578301504790424411472589676640}{
700438046366622594411215233016913071486089}}\,{T}^{6}$

\noindent$-{\frac {
193249058167562967816719266780463563218080}{
18911827251898810049102811291456652930124403}}\,{T}^{6}M+{\frac {
49806639391130680810149977028256722235840}{
1531858007403803613977327714607988887340076643}}\,{T}^{6}{M}^{3}$

\noindent$-{
\frac {6528982704325582783377044452994315443840}{
170206445267089290441925301623109876371119627}}\,{T}^{6}{M}^{2}+{
\frac {5956207676950329133758464677648797191680}{
1531858007403803613977327714607988887340076643}}\,{T}^{6}{M}^{4}$

\noindent$-{
\frac {687325278004245713061424259642842230089632}{
510619335801267871325775904869329629113358881}}\,{T}^{4}{M}^{2}-{
\frac {2050583245437666120312007255881326007529568}{
4595574022211410841931983143823966662020229929}}\,{T}^{4}{M}^{3}$

\noindent$+{
\frac {9808758961750027034269880327684876902994056}{
56735481755696430147308433874369958790373209}}\,{T}^{4}M-{\frac {
242134670257750264203544183229519465935616}{
4595574022211410841931983143823966662020229929}}\,{T}^{4}{M}^{4}$

\noindent$-{
\frac {396343814205392842122901681054670066164296}{
700438046366622594411215233016913071486089}}\,{T}^{4}+{\frac {
2035287795094287693870135081243409761340961}{
1531858007403803613977327714607988887340076643}}\,{T}^{2}{M}^{3}$

\noindent$+{
\frac {123744077281693368649685336016663155804936}{
1531858007403803613977327714607988887340076643}}\,{T}^{2}{M}^{4}+{
\frac {2168606272638745909274782927486819511187308}{
56735481755696430147308433874369958790373209}}\,{T}^{2}{M}^{2}$

\noindent$+{\frac 
{14909180233466477515674693647959213932639919}{
4202628278199735566467291398101478428916534}}\,{T}^{2}-{\frac {
33628313797235830289730557740669826423944793}{
37823654503797620098205622582913305860248806}}\,{T}^{2}M
 $.}
   $$ ~$$ 

  \begin{equation}\label{QUARUN51}
\left.\begin{aligned} \dot x &=
-y+a_{{1,1}}xy+{x}^{2}-a_{{1,1}}{x}^{3}
\\                  \dot y &=x+4\,a_{{1,1}}{y}^{2}-2\,xy-\frac{3}{2}\,a_{{1,1}}{x}^{2}-6\,a_{{1,1}}{x}^{2}y
+ \left( 2+{a_{{1,1}}}^{2} \right) {x}^{3}+ \left( 2\,a_{{1,1}}-\frac{1}{4}\,{
a_{{1,1}}}^{3} \right) {x}^{4}
 \end{aligned}\right\}, 
   \end{equation}
   $$ ~$$ 
   
  \begin{equation}\label{QUARUN52}
\left.\begin{aligned} \dot x &=-y+a_{{1,1}}xy+{x}^{2}-a_{{1,1}}{x}^{3}
\\                  \dot y &=x+3\,a_{{1,1}}{y}^{2}-2\,xy-a_{{1,1}}{x}^{2}-4\,a_{{1,1}}{x}^{2}y+
 \left( \frac{1}{3}\,{a_{{1,1}}}^{2}+2 \right) {x}^{3}+a_{{1,1}}{x}^{4} \end{aligned}\right\}, 
   \end{equation}
   $$ ~$$ 
   
  \begin{equation}\label{QUARUN53}
\left.\begin{aligned} \dot x &=-y+2\,\beta\,xy+{x}^{2}-2\,\beta\,{x}^{3}
\\                  \dot y &=x+8\,\beta\,{y}^{2}-2\,xy-3\,\beta\,{x}^{2}-12\,\beta\,{x}^{2}y+14\,{x}^{3}-2\,\beta\,{x}^{4}
 \end{aligned}\right\}, 
   \end{equation}
\noindent where $\beta=\pm \,\sqrt{3}$.
  $$ ~$$ 
  
  \begin{equation}\label{QUARUN54}
\left.\begin{aligned} \dot x &=-y+\alpha\,xy+{x}^{2}-4/3\,\alpha\,{x}^{3}+2/3\,{x}^{4}
\\                  \dot y &=x+6\,\alpha\,{y}^{2}-2\,xy-5/2\,\alpha\,{x}^{2}-9\,\alpha\,{x}^{2}y+{
\frac {26}{3}}\,{x}^{3}+6\,{x}^{3}y-8/3\,\alpha\,{x}^{4}\end{aligned}\right\},  
   \end{equation}
   $$ ~$$ 
   
  \begin{equation}\label{QUARUN55}
\left.\begin{aligned} \dot x &=-y+\alpha\,xy+{x}^{2}-\frac{4}{3}\,\alpha\,{x}^{3}+\frac{2}{3}\,{x}^{4}
\\                  \dot y &=x+3\,\alpha\,{y}^{2}-2\,xy-\alpha\,{x}^{2}-3\,\alpha\,{x}^{2}y+\frac{8}{3}\,{x
}^{3}+2\,{x}^{3}y-\frac{2}{3}\,\alpha\,{x}^{4} \end{aligned}\right\}, 
\end{equation}
\noindent where $\alpha=\pm\,\sqrt{2}$.
  $$ ~$$ 
  
  \begin{equation}\label{QUARUN56}
\left.\begin{aligned} \dot x &=-y+a_{{1,1}}xy+{x}^{2}+ \left( {a_{{1,1}}}^{2}-\frac{3}{2}\,b_{{0,2}}a_{{1,1}}
+\frac{1}{2}\,{b_{{0,2}}}^{2} \right) {x}^{2}y-a_{{1,1}}{x}^{3}+ \left( -{a_{{
1,1}}}^{2}+\frac{3}{2}\,b_{{0,2}}a_{{1,1}}-\frac{1}{2}\,{b_{{0,2}}}^{2} \right) {x}^{4
}
\\                  \dot y &=x+b_{{0,2}}{y}^{2}-2\,xy+ \left( \frac{1}{2}\,a_{{1,1}}-\frac{1}{2}\,b_{{0,2}}
 \right) {x}^{2}+ \left( 2\,{a_{{1,1}}}^{2}-3\,b_{{0,2}}a_{{1,1}}+{b_{
{0,2}}}^{2} \right) x{y}^{2}\\
&+ \left( 2\,a_{{1,1}}-2\,b_{{0,2}}
 \right) {x}^{2}y+2\,{x}^{3}+ \left( -2\,{a_{{1,1}}}^{2}+3\,b_{{0,2}}a
_{{1,1}}-{b_{{0,2}}}^{2} \right) {x}^{3}y+ \left( -2\,a_{{1,1}}+b_{{0,
2}} \right) {x}^{4}
 \end{aligned}\right\}, 
   \end{equation}
   $$ ~$$ 
   
    \begin{equation}\label{QUARUN57}
\left.\begin{split} \dot x &=-y+\alpha\,xy+M\,{x}^{2}y+{x}^{2}-\alpha\,{x}^{3}-M\,{x}^{4}\\
\dot y&=x+5\,Mx{y}^{2}-8\,M{x}^{3}y -2\,xy+\frac{B_{{0,2}}}{12\,({ \left( M+9 \right) \alpha\,\beta})} {y}^{2}-{\frac {\delta}{ 12\left( M+9
 \right) \alpha\,\beta}} {x}^{2}\\
 &-\,{\frac {\delta}{ 3\left( M+9 \right) \alpha\,\beta}}{x}^{2}y-\frac{B_{{3,0}}}{2 \left( M+9 \right) \beta}{x}^{3}- \frac{B_{{4,0}}}{12{
 \left( M+9 \right) \alpha\,\beta}} {x}^{4}\end{split}\right\}, 
   \end{equation}

\noindent where $\left( M+9 \right) \alpha\,\beta\neq 0$ and $\alpha$ is the real root of the polynomial $P(Z):$ 
  $$ ~$$ 
\noindent$P(Z)= \left( 144\,M+1296 \right) {{
\it \_Z}}^{6}+ \left( 7236\,M-24\,{M}^{3}+996\,{M}^{2}-15552 \right) {{\it \_Z}}^{4}$

\noindent$+ \left( 12960\,{M}^{2}+
2736\,{M}^{3}-202\,{M}^{4}-46656\,M\right) {{
\it \_Z}}^{2}$

\noindent$-8\,{M}^{6}+664\,{M}^{4}-109\,{M}
^{5}-15552\,{M}^{2}-233280\,M+10368\,{M}^{3}$
\\
and \\
$\beta=1581568416-224710848\,{M}^{3}+452213280\,{M}^{2}-
20366712\,{M}^{4}+284739\,{M}^{6}+3438793\,{M}^{5}+32\,{M}^{9}-892\,{M
}^{8}-15958\,{M}^{7}+2301021432\,{\alpha}^{2}+655618860\,{\alpha}^{4}+
4745755008\,M+96\,{\alpha}^{4}{M}^{6}-40132\,{\alpha}^{2}{M}^{6}+1000
\,{M}^{7}{\alpha}^{2}-3120\,{\alpha}^{4}{M}^{5}+128736\,{\alpha}^{4}{M
}^{4}+530942\,{\alpha}^{2}{M}^{5}-178814952\,{M}^{2}{\alpha}^{2}+
4111938\,{M}^{4}{\alpha}^{2}+1440585648\,{\alpha}^{2}M-68912208\,{M}^{
3}{\alpha}^{2}-36318780\,M{\alpha}^{4}+216852\,{M}^{3}{\alpha}^{4}-
22879908\,{M}^{2}{\alpha}^{4}$,
$$ ~$$
$\delta=-1309805864856\,{M}^{3}+
1304677734768\,{M}^{2}-109564995222\,{M}^{4}+3102777711\,{M}^{6}+
32534070309\,{M}^{5}+288136\,{M}^{9}-31013013\,{M}^{8}-281741535\,{M}^
{7}+44302\,{M}^{10}+341618777856\,{\alpha}^{2}+1640171477952\,{\alpha}
^{4}+1000\,{M}^{10}{\alpha}^{2}+96\,{\alpha}^{4}{M}^{9}-1680\,{\alpha}
^{4}{M}^{8}+17574286201920\,M+1910484\,{\alpha}^{2}{M}^{8}-49132\,{
\alpha}^{2}{M}^{9}+4176\,{\alpha}^{4}{M}^{7}-4906584\,{\alpha}^{4}{M}^
{6}-286521084\,{\alpha}^{2}{M}^{6}-30232362\,{M}^{7}{\alpha}^{2}+
34671168\,{\alpha}^{4}{M}^{5}+1921156704\,{\alpha}^{4}{M}^{4}+
7550896788\,{\alpha}^{2}{M}^{5}-708776248212\,{M}^{2}{\alpha}^{2}+
32149733976\,{M}^{4}{\alpha}^{2}+4446970742160\,{\alpha}^{2}M-
333064034298\,{M}^{3}{\alpha}^{2}-226136884752\,M{\alpha}^{4}+
4800914064\,{M}^{3}{\alpha}^{4}-99633193992\,{M}^{2}{\alpha}^{4}-796\,
{M}^{11}+32\,{M}^{12}$,

$$ ~$$

$ B_{{0,2}}=
854046944640\,{\alpha}^{2}+4378535313120\,{\alpha}^{4}+18416851416\,{M}^{3}{\alpha}^{4}-263459876832\,{M}^{2}{\alpha}^{4}-
711485573796\,M{\alpha}^{4}+11974507973856\,{\alpha}^{2}M-884184180804\,{M}^{3}{\alpha}^{2}+104273686746\,{M
}^{4}{\alpha}^{2}-2161040187528\,{M}^{2}{\alpha}^{2}+
20100360978\,{\alpha}^{2}{M}^{5}+5321931120\,{\alpha}^{4}{M}^{4}+30597468\,{\alpha}^{4}{M}^{5}-16021728\,{\alpha}^{4
}{M}^{6}-1086777654\,{\alpha}^{2}{M}^{6}-80729938\,{M}^{7}{\alpha}^{2}+93312\,{\alpha}^{4}{M}^{7}-179920\,{
\alpha}^{2}{M}^{9}+6307872\,{\alpha}^{2}{M}^{8}-5568\,
{\alpha}^{4}{M}^{8}+384\,{\alpha}^{4}{M}^{9}+4000\,{M}^{10}{\alpha}^{2}-3677920975152\,{M}^{3}-249464899332
\,{M}^{4}+7497589581\,{M}^{6}+95584892967\,{M}
^{5}+140776\,{M}^{10}+1979888\,{M}^{9}-79612855\,{M}^{8}-920925749\,{M}^{7}-2800\,{M}^{11}+128\,{M}^{12}+47893803042240\,M+2753000428896\,{M}^{2}$,
$$ ~$$
$B_{{3,0}}=-
56936462976-139773234528\,{\alpha}^{2}-204261985872\,{\alpha}^{4}+
14956272\,{M}^{3}{\alpha}^{4}+11449948104\,{M}^{2}{
\alpha}^{4}+13633897968\,M{\alpha}^{4}-489497786712\,{\alpha}^
{2}M+35663387130\,{M}^{3}{\alpha}^{2}-1670720580\,{M}^{4}{\alpha}^{2}+51974431164\,{M}^{2}{\alpha}^{2}-
718427376\,{\alpha}^{2}{M}^{5}-165970176\,{\alpha}^{4}{M}^{4}-5669352\,{\alpha}^{4}{M}^{5}+260880\,{\alpha}^{4}{M}^{6}+8347234\,{\alpha}^{2}{M}^{6}+2613748\,{M}^{
7}{\alpha}^{2}-3696\,{\alpha}^{4}{M}^{7}+3000\,{\alpha}^{2}{M}^{9}-119572\,{\alpha}^{2}{M}^{8}+288\,{\alpha}^{4}{M}^{8}+117311420664\,{M}^{3}+14135146866\,{M}^{4}
-342595189\,{M}^{6}-2512013814\,{M}^{5}-1940\,{M}^{10}+36898\,{M}^{9}+2916172\,{M}^{8}+15602085\,{M}^{7}+96\,{M}^{11}-1772624526048\,M-
200707918416\,{M}^{2},$
$$ ~$$
$B_{{4,0}}=-512428166784\,{\alpha}^
{2}+1155128820192\,{\alpha}^{4}+76187990232\,{M}^{3}{\alpha}^{
4}-37449803088\,{M}^{2}{\alpha}^{4}-1560658377780\,M{
\alpha}^{4}+4471598426688\,{\alpha}^{2}M-170217184320\,{M}^{3}{\alpha}^{2}+262466972514\,{M}^{4}{\alpha}^{2}-
3995129998656\,{M}^{2}{\alpha}^{2}+4574201922\,{\alpha}^{2}{M}^{5}+3865776624\,{\alpha}^{4}{M}^{4}-781052628\,{
\alpha}^{4}{M}^{5}-41458608\,{\alpha}^{4}{M}^{6}-
4386392646\,{\alpha}^{2}{M}^{6}-21588886\,{M}^{7}{
\alpha}^{2}+1071072\,{\alpha}^{4}{M}^{7}-668472\,{\alpha}^{2}{
M}^{9}+17045880\,{\alpha}^{2}{M}^{8}-15264\,{\alpha}^{
4}{M}^{8}+1728\,{\alpha}^{4}{M}^{9}+18000\,{M}
^{10}{\alpha}^{2}-3279573271872\,{M}^{3}+482166146232\,{M}^{4}-8025777621\,{M}^{6}+136445218065\,{M}^{5}+
323820\,{M}^{10}+15941920\,{M}^{9}+19475327\,{M}^{8}-2392822547\,{M}^{7}-9336\,{M}^{11}+576\,{M}^{12}+25093708683840\,M-8570037406464\,{M}^{2}.$

\end{thm}

\section{Explicit Linearization}\label{sec::lin}
\subsection{Linearization formulas}
 Let us consider the Liénard type system \eqref{LienardType}
 \begin{equation*}
\left.\begin{aligned} \dot x &= y \\ \dot y &= -g(x) - f(x) y^2
 \end{aligned}\right\} 
   \end{equation*} 
  with a center at the origin $(0,0)$ where $f$ and $g$ are real analytic in a neighborhood of zero.
  

It is known by Sabatini formula \eqref{L2FI} that the first integral associated to the system~\eqref{LienardType} can be written
\begin{equation}\label{SAB}
I(x,\dot x)= \int_0^x g(s) e^{2F(s)} ds +\frac{1}{2}(\dot x e^{F(x)})^2
\end{equation}
where $F(x)=\int_0^x f(s) ds$. 

Following \cite{ChoudhuryGuha2010} (see also \cite{GuhaChoudhury2009}), let us perform the following change of variables
  \begin{equation}\label{CNG}
\left.\begin{aligned} p(x,\dot x)&=\dot x e^{F(x)} \\ q(x)&=\int_0^x e^{F(s)} ds
 \end{aligned}\right\}. 
   \end{equation}

   As $\dfrac{\partial(p,q)}{\partial(x,\dot x)}=-e^{2F(x)}<0$ and $p(0,0)=q(0)=0$ then this is a change of variables preserving the origin and well defined around it. Moreover, $q'(x)=e^{F(x)}>0$ and thus the function $x\longmapsto q(x)$ is strictly increasing. In the $(p,q)$ coordinates the first integral~\eqref{SAB} becomes
  \begin{equation}\label{HAM}I(x,\dot x)=H(p,q)=\frac{1}{2}p^2+U(q),
  \end{equation}
  where $U$ is some uniquely defined real analytic function, $U(0)=0$.
  Now it is easy to see that the system \eqref{LienardType} in $(p,q)$ coordinates can be written as
  
  \begin{equation}\label{POT}
\left.\begin{aligned} \dot q&=\dfrac{\partial H}{\partial p}= p \\ \dot p &=-\dfrac{\partial H}{\partial q}=- \dfrac{d }{d q}U
 \end{aligned}\right\} 
  \end{equation} 
  that is as a Hamiltonian system corresponding to the Hamiltonian~\eqref{HAM}.
  
  The main result of this Section is 
  
  \begin{thm}\label{thm:Linea}
  Let us consider the Li\'enard type system~\eqref{LienardType} with real analytic functions $f$ and $g$ such that $x\,g(x)>0$ for $x\neq 0$. Then the origin $O$ is an isochronous center with Urabe function $h=0$ if and only if $U(q)=\frac{q^2}{2}$. 
\end{thm}  
\begin{proof}
Taking into account our assumptions, from Theorem \ref{thm:Center} one knows that $O$ is a center. Now, from Corollary \ref{cor:urabenul}, one knows that $O$ is an isochronous center with Urabe function $h=0$ if and only if  $ g'(x)+g(x)f(x)=1$ or equivalently $g'(x)e^{F(x)}+g(x)f(x)e^{F(x)}=e^{F(x)}$. The last equality is nothing else $\left(g(x)e^{F(x)}\right)'=\left(\int_0^x e^{F(s)}ds\right)'$. As $g(0)e^{F(0)}=0$, when integrating one obtains $g(x)e^{F(x)}=\int_0^x e^{F(s)}ds$ or equivalently $U'(q)=q$ because $(\frac{dU}{dq})(q(x))=g(x)\,e^{F(x)}$. Since $U(0)=0$ one has $U(q)=\frac{1}{2}q^2$.
\end{proof}
 Consequently when the Urabe function $h$ identically vanishes , the system of coordinate $(p,q)$ defined by \eqref{CNG} is the linearizing system of coordinates for system~\eqref{LienardType}. Indeed, 
\begin{equation}\label{LIN}
\left.\begin{aligned} \dot q&= p \\ \dot p &=-q
 \end{aligned}\right\}. 
  \end{equation} 
  
  \medskip
 It is interesting to compare the theorem \ref{thm:Linea} with Chalykh-Veselov theorem  that we formulate only for potential $U$ without pole at $0$.
 
 \begin{thm}[\cite{ChalykhVeselov2005}, Theorem 1]
  Let us consider the Hamiltonian system with the Hamiltonian $H(P,Q)=\frac{1}{2}P^2+U(Q)$ where $U$ is a rational function without pole at $0$. Then $O$ is an isochronous center for the associated Hamiltonian system $\dot Q=\dfrac{\partial H}{\partial P},\; \dot P =-\dfrac{\partial H}{\partial Q}$ if and only if up to a shift $Q\to Q+a$ and adding a constant, $U(Q)=k\,Q^2$ for some $k\in \mathbb{R}-\{0\}$.
  \end{thm}  

  \subsection{Examples }
  Now applying formula~\eqref{CNG} we provide $5$ examples of linearization of isochronous centers with zero Urabe function. The reduction to Li\'enard type system \eqref{LienardType} is always obtained by standard or Choudhury-Guha reduction. To compute variables $(p,q)$ see (\eqref{CNG}) we use Maple, and the identity \eqref{LIN} was verified in all cases. Our choice is somewhere random, because all reported examples with zero Urabe function are good for such purpose.
  
  \subsubsection{Cubic examples} 
 \begin{enumerate}
 \item 
 Consider the case 6 of Theorem~3
  of~\cite{ChouikhaRomanovskiChen2007}, that is the system
\begin{equation}\label{CUBI}
\left.\begin{aligned} \dot x&=-y+a_{{1,1}}xy+ \left( {a_{{1,1}}}^{2}-\frac{3\,a_{{1,1}}b_{{0,2}}}{2}+\frac{{b_{{0,
2}}}^{2}}{2} \right) {x}^{2}y
\\ \dot y &=x+ \left( -\frac{b_{{0,2}}}{2}+\frac{a_{{1,1}}}{2} \right) {x}^{2}+b_{{0,2}}{y}^{
2}+ \left( 2\,{a_{{1,1}}}^{2}-3\,a_{{1,1}}b_{{0,2}}+{b_{{0,2}}}^{2}
 \right) x{y}^{2}
 \end{aligned}\right\}. 
  \end{equation}  
  In this case the functions $f$ and $g$ are 
  \begin{equation}\label{FGCUB}
\left.\begin{aligned} f(x)&=-{\frac {b_{{0,2}}+ \left( 2\,{a_{{1,1}}}^{2}-3\,a_{{1,1}}b_{{0,2}}+{b_{
{0,2}}}^{2} \right) x+a_{{1,1}}+2\, \left( {a_{{1,1}}}^{2}-\frac{3\,a_{{1,1}}b_{{0,2}}}{2}+\frac{{b_{{0,
2}}}^{2}}{2} \right) x}{-1+a_{{1,1}}x+ \left( {a_{{1,1
}}}^{2}-\frac{3\,a_{{1,1}}b_{{0,2}}}{2}+\frac{{b_{{0,
2}}}^{2}}{2} \right) {x}^{2}}}
\\ g(x) &= \left( 1-a_{{1,1}}x- \left( {a_{{1,1}}}^{2}-\frac{3\,a_{{1,1}}b_{{0,2}}}{2}+\frac{
{b_{{0,2}}}^{2}}{2} \right) {x}^{2} \right)  \left( x+ \left( -\frac{b_{{0,2}}}{2}+\frac{a_{{1,1}}}{2} \right) {x}^{2} \right) 
 \end{aligned}\right\}, 
  \end{equation}  
  for which we obtain the following linearizing change of coordinates
  \begin{equation}
\left.\begin{aligned} q(x)&=- \dfrac{\left( 2+a_{{1,1}}x-b_{{0,2}}x \right) x{{\rm e}^{2\, \left(a_{{1,1}} -b_{{0
,2}} \right)  \frac{\left( A(x) -A(0)  \right)} {{\sqrt {-5\,{a_{{
1,1}}}^{2}+6\,a_{{1,1}}b_{{0,2}}-2\,{b_{{0,2}}}^{2}}}}}}}{ \left( -2+2\,a_
{{1,1}}x+2\,{a_{{1,1}}}^{2}{x}^{2}-3\,{x}^{2}a_{{1,1}}b_{{0,2}}+{x}^{2}{b
_{{0,2}}}^{2} \right)}\\
p(x,y) &=\dfrac{-2y\,q(x)}{\left( 2+a_{{1,1}}x-b_{{0,2}}x \right)x}
 \end{aligned}\right\} 
  \end{equation}  
  where
  \begin{equation}
A(x)=\arctan \left( {\frac {2\,{a_{{1,1}}}^{2}
x-3\,x\,a_{{1,1}}b_{{0,2}}+x{b_{{0,2}}}^{2}+a_{{1,1}}}{\sqrt {-5\,{a_{{1,1}
}}^{2}+6\,a_{{1,1}}b_{{0,2}}-2\,{b_{{0,2}}}^{2}}}} \right).
  \end{equation}  
  
  \item Consider system \eqref{CUB1} from Theorem~\ref{CUBICS}. 
  In this case the functions $f$ and $g$ are 
  \begin{equation}\label{FGCUBII}
\left.\begin{aligned} f(x)&=-6\,{\frac {b_{{2,0}}}{1+2\,b_{{2,0}}x}}\\ 
g(x) &= x \left( 2\,{b_{{2,0}}}^{2}{x}^{2}+3\,b_{{2,0}}x+1 \right)
 \end{aligned}\right\}, 
  \end{equation}  
  for which we obtain the following linearizing change of coordinates
  \begin{equation}
\left.\begin{aligned} q(x)&={\frac {x \left( b_{{2,0}}x+1 \right) }{ \left( 1+2\,b_{{2,0}}x
 \right) ^{2}}}
\\
p(x,y) &={\frac {-y+{x}^{2}}{ \left( 1+2\,b_{{2,0}}x \right) ^{2}}}
 \end{aligned}\right\}. 
  \end{equation}  

  \end{enumerate} 
  
   \subsubsection{Quartic examples}
  \begin{enumerate}
  \item As a quartic example we consider the system III of Theorem 3.1 of~\cite{BoussaadaChouikhaStrelcyn2010}. We choose the following restrictions on the parameters $ a_{{21}}=-3,b_{{0,2}}=-4,b_{{20}}=\frac{1}{2} $ to obtain simple,  presentable expressions for linearizing variables.
   \begin{equation}\label{QUA}
\left.\begin{aligned} \dot x&=-y-3\,xy-3\,{x}^{2}y-{x}^{3}y
\\ \dot y &=x+\frac{1}{2}\,{x}^{2}-4\,{y}^{2}-4\,x{y}^{2}-2\,{x}^{2}{y}^{2}
 \end{aligned}\right\}. 
  \end{equation}  
  In this case the functions $f$ and $g$ are 
  \begin{equation}
\left.\begin{aligned} f(x)&=-{\frac {7+10\,x+5\,{x}^{2}}{1+3\,x+3\,{x}^{2}+{x}^{3}}}
\\ g(x) &= \frac{1}{2}\, \left( 1+3\,x+3\,{x}^{2}+{x}^{3} \right) x \left( 2+x \right) 
 \end{aligned}\right\}, 
  \end{equation}  
  for which we obtain the following linearizing change of coordinates
  \begin{equation}
\left.\begin{aligned} q(x)&=\frac{1}{2}\,x \left( 2+x \right) {{\rm e}^{-{\frac {x \left( 2+x \right) }{
 \left( 1+x \right) ^{2}}}}} \left( 1+x \right) ^{-2}
\\
p(x, y)&= -y{{\rm e}^{-{\frac {x \left( 2+x \right) }{ \left( 1+x \right) ^{2}}}
}} \left( 1+x \right) ^{-2}
 \end{aligned}\right\}. 
  \end{equation}  
  
   \item  We consider the system \eqref{QUARUN42} of Theorem \ref{thm:QUARTICURABENUL}.

  In this case the functions $f$ and $g$ are 
  \begin{equation}
\left.\begin{aligned} f(x)&={\frac {4\,{a_{{1,1}}}^{2}x-6\,xb_{{0,2}}a_{{1,1}}+2\,x{b_{{0,2}}}^{2}
+a_{{1,1}}+b_{{0,2}}}{1-a_{{1,1}}x-{a_{{1,1}}}^{2}{x}^{2}+3/2\,{x}^{2}
b_{{0,2}}a_{{1,1}}-1/2\,{x}^{2}{b_{{0,2}}}^{2}}}\\
g(x)&=-1/4\,x \left( a_{{1,1}}x-xb_{{0,2}}+2 \right)  \left( -2+2\,a_{{1,1}}
x+2\,{a_{{1,1}}}^{2}{x}^{2}-3\,{x}^{2}b_{{0,2}}a_{{1,1}}+{x}^{2}{b_{{0
,2}}}^{2} \right) 
 \end{aligned}\right\}, 
  \end{equation}  
  for which we obtain the following linearizing change of coordinates
  \begin{equation}
\left.\begin{aligned} q(x)&=-{\frac {x \left( a_{{1,1}}x-xb_{{0,2}}+2 \right) S \left( x \right) }
{-2+2\,a_{{1,1}}x+2\,{a_{{1,1}}}^{2}{x}^{2}-3\,{x}^{2}b_{{0,2}}a_{{1,1
}}+{x}^{2}{b_{{0,2}}}^{2}}}
\\
p(x, y)&= 2\,{\frac { \left( {x}^{2}-y \right) S \left( x \right) }{-2+2\,a_{{1,
1}}x+2\,{a_{{1,1}}}^{2}{x}^{2}-3\,{x}^{2}b_{{0,2}}a_{{1,1}}+{x}^{2}{b_
{{0,2}}}^{2}}}
 \end{aligned}\right\}, 
  \end{equation}  
  \end{enumerate}
  
  where $$ S \left( x \right) ={{e}^{{\frac {2\, \left( a_{{1,1}}-b_{{0,2}} \right) 
 \left( \arctan \left( {\frac {2\,{a_{{1,1}}}^{2}x-3\,xb_{{0,2}}a_{{1,
1}}+x{b_{{0,2}}}^{2}+a_{{1,1}}}{\sqrt {-5\,{a_{{1,1}}}^{2}+6\,b_{{0,2}
}a_{{1,1}}-2\,{b_{{0,2}}}^{2}}}} \right) -\arctan \left( {\frac {a_{{1
,1}}}{\sqrt {-5\,{a_{{1,1}}}^{2}+6\,b_{{0,2}}a_{{1,1}}-2\,{b_{{0,2}}}^
{2}}}} \right)  \right) }{\sqrt {-5\,{a_{{1,1}}}^{2}+6\,b_{{0
,2}}a_{{1,1}}-2\,{b_{{0,2}}}^{2}}}}}}.$$

   \subsubsection{A rational example}
    Let us consider the system
    \begin{equation}\label{RAT}
\left.\begin{aligned} \dot x&=-y+ \dfrac{yx}{1+x}
\\ \dot y &=x+\dfrac{y^2}{1+x}
 \end{aligned}\right\}.
  \end{equation}  
  In this case the functions $f$ and $g$ are
  \begin{equation}
\left.\begin{aligned} f(x)&={\dfrac {2+x}{1+x}}
\\ g(x) &= {\dfrac {x}{1+x}} 
 \end{aligned}\right\}, 
  \end{equation}  
  for which we obtain the following linearizing change of coordinates
  \begin{equation}
\left.\begin{aligned} q(x)&=x{{\rm e}^{x}}
\\
p(x, y)&=y{{\rm e}^{x}} \left( 1+x \right)
 \end{aligned}\right\}.
  \end{equation}  
\subsection{Comments}
It is really astonishing that in all above cases the linearizing variables $(p,q)$ are always expressed in ''finite terms''. This follows from the fact that if  $g'(x)+f(x)\,g(x)=1$ then $f=\frac{1-g'}{g}$. Moreover, as $g(0)=0$ and $g'(0)=1$ the singularity of $f$ at zero is spurious. In all examples considered in this and related papers \cite{Chouikha2007,ChouikhaRomanovskiChen2007,BoussaadaChouikhaStrelcyn2010} $f$ and $g$ always are  rational functions. Then $F(x)=\int_0^x f(s)\,ds$ is expressed in ''finite terms'' and thus also $p(x,\dot x)$. The problem is slightly more delicate in what concerns $q(x)$. But $\int e^{F(s)}\,ds=\int e^{\int f(s)\,ds}\,ds=\int e^{\int\frac{1-g'(s)}{g(s)}\,ds}\,ds=\int\frac{ 1}{g(s)}e^{\int \frac{1}{g(s)}\,ds}\,ds=e^{\int \frac{1}{g(s)}\,ds}+const$. $g$ being a rational function, $\int \frac{ds}{g(s)}$ is obtained in ''finite terms'' and thus also $\int e^{F(s)}\,ds$.

\subsection*{Acknowledgment}
We are grateful to the CRIHAN (Centre de Ressources Informatiques de
Haute Normandie), which provided us technical resources for our
computations.

We also sincerely thank Marie-Claude Werquin (University Paris 13, France)
who corrected and improved our scientific English. 
\bibliographystyle{plain}
\bibliography{BBCS}

\end{document}